\documentclass[preprint,11pt]{article}
\pdfoutput=1
\usepackage{amssymb,amsfonts,amsmath,amsthm,amscd,dsfont,mathrsfs}
\usepackage{graphicx,float,psfrag,epsfig}
\usepackage{wrapfig}
\usepackage{relsize}
\usepackage{color}
\usepackage{pict2e}
\usepackage[tight]{subfigure}
\usepackage{algorithm}
\usepackage[noend]{algorithmic}
\usepackage{caption}

\DeclareMathAlphabet{\mathpzc}{OT1}{pzc}{m}{it}

\newtheorem{propo}{Proposition}[section]
\newtheorem{lemma}[propo]{Lemma}
\newtheorem{definition}[propo]{Definition}
\newtheorem{coro}[propo]{Corollary}
\newtheorem{thm}[propo]{Theorem}

\def\cF{{\cal F}}
\def\cA{{\cal A}}
\def\cB{{\cal B}}
\def\cR{{\cal R}}

\def\cT{{\cal T}}

\def\reals{{\mathbb R}}

\def\dual{\Psi}

\def\eps{{\varepsilon}}

\def\E{{\mathbb E}}
\def\Var{{\rm Var}}

\def\L0{{L_0}}

\def\de{{\rm d}}
\def\<{\langle}
\def\>{\rangle}

\def\bX{{\mathbf X}}

\def\htheta{\widehat{\theta}}

\def\limsup{{\rm lim \, sup}}

\def\F{{\sf F}}

\def \Tr{{\rm Trace}}
\def\F{{\sf F}}

\def\P{{\mathbb{P}}}

\def\sT{{\sf T}}

\def\v*{v_0}
\def\T*{T_0}

\def\u*{u_0}
\def\F*{F_0}

\definecolor{olivegreen}{rgb}{0,0.6,0.4}

\def\AMSE{{\rm AMSE}}
\def\dof{{m}}
\def\tkappa{\tilde{\kappa}}

\newcommand{\bitem}{\begin{itemize}}
\newcommand{\eitem}{\end{itemize}}
\newcommand{\goto}{\to}
\newcommand{\beq}{\begin{equation}}
\newcommand{\eeq}{\end{equation}}

\def\tarr{\bar{r}}
\def\barr{{\bar{r}}}
\def\bbarr{\bar{\bar{r}}}
\def\bbtau{\bar{\bar{\tau}}}
\def\bkappa{\bar{\bar{\kappa}}}
\def\bbcA{\bar{\bar{\cA}}}
\def\bcT{\bar{\cT}}
\def\bbcT{\bar{\bcT}}
\def\ukappa{\underline{\kappa}}

\def\cV{{\cal V}}
\def\cK{{\cal K}}
\def\Var{{\rm Var}}
\def\AVar{{\rm AVar}}
\def\Tr{{\rm Tr}}

\newcommand{\ajcomment}[1]{}

\makeatletter
\newcommand{\labitem}[2]{%
\def\@itemlabel{\text{#1}}
\item
\def\@currentlabel{#1}\label{#2}}

\usepackage{bibentry}
\newcommand{\ignore}[1]{}
\newcommand{\nobibentry}[1]{{\let\nocite\ignore\bibentry{#1}}}

\makeatother

\addtocontents{toc}{\protect\setcounter{tocdepth}{2}}

\title{Variance Breakdown of Huber $(M)$-estimators: $n/p \goto m \in (1,\infty)$}

\author{David Donoho
            \footnote{Department of Statistics, Stanford
              University} \, and \,Andrea~Montanari 
            \footnote{Department of Electrical Engineering and
              Department of Statistics, Stanford University}
            }

\begin{document}

\date{March 5, 2015}

\maketitle
%\thispagestyle{empty}
%\pagestyle{empty}

%%%%%%%%%%%%%%%%%%%%%%%%%%%%%%%%%%%%%%%%%%%%%%%%%%%%%%%%%%%%%%%%%%%%%%%%%%%%%%%%
\begin{abstract}
Huber's gross-errors contamination model 
considers the class $\cF_\eps$ of all noise distributions
$F = (1 - \eps) \Phi + \eps H$, with $\Phi$ standard normal,
$\eps \in (0,1)$ the contamination fraction,
and $H$ the contaminating distribution. A half century ago, Huber
evaluated the minimax asymptotic variance in scalar location estimation,
\beq \label{eq:abstractClassical}
     \min_\psi \max_{F \in \cF_\eps}  V(\psi, F) = \frac{1}{I(F_\eps^*)}
\eeq
where $V(\psi,F)$ denotes the asymptotic variance of the $(M)$-estimator
for location with score function $\psi$, 
and $I(F_\eps^*)$ is the minimal Fisher information $ \min_{\cF_\eps} I(F)$.

We consider the linear regression model  $Y = X\theta_0 + W$, $W_i\sim_{\text{i.i.d.}}F$,
and iid Normal predictors $X_{i,j}$, working
in the high-dimensional-limit asymptotic  where the number $n$ of observations
and $p$ of variables   both grow large, while
$n/p \goto \dof \in (1,\infty)$; hence $\dof$ plays the role of `asymptotic number
of observations per parameter estimated'.  
%Classical formulas for asymptotic variance of $(M)$-estimates are
%incorrect in this setting;  Namely if
%$V(\psi,F)$ denotes the classical asymptotic variance in location estimation and  
Let $V_\dof(\psi,F)$ denote the per-coordinate asymptotic variance of the $(M)$-estimator of regression
in the $n/p \goto \dof$ regime \cite{karoui2013robust,donoho2013high,karoui2013cavity}.
Then $V_\dof \neq V$; however $V_\dof \goto V$ 
as $\dof \goto \infty$.

In this paper we evaluate the minimax asymptotic variance  of the Huber $(M)$-estimate.
The statistician minimizes over
the family $(\psi_\lambda)_{\lambda > 0}$ of all tunings of Huber $(M)$-estimates of regression,
and Nature maximizes over gross-error contaminations  $F \in \cF_\eps$. Suppose that $I(F_\eps^*) \cdot \dof > 1$.
Then
\beq \label{eq:abstractModern}
    \min_\lambda \max_{F \in \cF_\eps}  V_\dof(\psi_\lambda, F) =  \frac{1}{I(F_\eps^*) - 1/\dof} .
\eeq
Of course, the RHS of (\ref{eq:abstractModern})   is strictly bigger than the RHS of  (\ref{eq:abstractClassical}).
Strikingly, if  $I(F_\eps^*) \cdot \dof \leq  1$, then
\[
\min_\lambda \max_{F \in \cF_\eps}  V_\dof(\psi_\lambda, F) =  \infty .
\]
In short, the asymptotic variance of the Huber estimator
{\em breaks down} at a critical ratio of observations per parameter.
Classically, for the minimax $(M)$-estimator of  location, 
no such breakdown occurs \cite{donoho1983notion}. However, under this paper's
$n/p \goto m$ asymptotic,  the
breakdown point is where the Fisher information per  parameter equals unity:
\[
    \eps^* \equiv \eps_\dof^*(\mbox{Minimax Huber-} (M) \mbox{ Estimate}) =  \inf \{ \eps : \dof \cdot I(F_\eps^*) \geq 1 \}.
\]

\end{abstract}

{{\bf Dedication.} \sl Based on a lecture delivered at a special colloquium honoring
the 50th anniversary of the Seminar f\"ur Statistik (SfS) at ETH Z\"urich, November 25, 2014.
The year 2014 was simultaneously: the 80th birthday year of Peter Huber, the 50th anniversary of his great 1964
paper on Robust Estimation, and the 50th anniversary of SfS.  All of these events are causes
for celebration, and we thank especially Peter B\"uhlmann,  Sara van de Geer,
Hansruedi K\"unsch,  Marloes Maathuis, Nicolai Meinshausen,
and indeed everyone at SfS for creating a wonderful commemoration
event. Special congratulations to Peter J. Bickel on receiving his Doctor Honoris
Causa from ETH as part of this celebration!}

%%%%%%%%%%%%%%%%%%%%%%%%%%%%%%%%%%%%%%%%%%%%%%%%%%%%%%%%%%%%%%%%%%%%%%%%%%%%%%%%
\section{Introduction}

Fifty years ago, Peter Huber published the masterwork  \cite{HuberMinimax}
in the {\em Annals of Mathematical Statistics}.
His paper, `Robust Estimation of a Location Parameter'
revealed robust statistics to be amenable to mathematical analysis,
producing a new optimal   robust estimator --
now called the Huber (M)-estimator -- that has proven 
 practical, elegant and lasting.  Richard Olshen
once called Peter's paper `an out-of-the-park, grand-slam
home run'.\footnote{Terminology from American
baseball. The highest-impact scoring outcome that can ever be delivered by a
batsman, and not at all frequent. Wikipedia states that over 112 annual World Series,
comprising more than 500 games, and ten thousand at-bats,
this has happened only
eighteen times. }

Only 8 years after this initial paper in statistics,
Peter delivered the Wald Lectures \cite{huber1973robust},  
recognition from the profession of the exceptional importance
of his {\em \oe uvre}. While Huber's 1964 paper
considered the estimation of a scalar location parameter,   
his Wald Lectures summarized work showing that
much of the framework of the  1964
paper  generalized immediately to
regression estimation.

\subsection{(M)-estimates of Regression}

Consider the traditional  linear regression model
\begin{eqnarray}\label{eq:NoisyModel}
Y\, =\, \bX\,\theta_0+ W\, ,\label{eq:Model}
\end{eqnarray}
with  $Y = (Y_1,\dotsc,Y_n)^\sT\in\reals^n$ a vector of
responses, $\bX\in \reals^{n\times p}$  a known design matrix,
$\theta_0\in\reals^p$ a vector of parameters, and $W\in\reals^n$ 
a random noise vector with i.i.d. components 
 having marginal  distribution $F = F_{W}$.\footnote{With a slight abuse of notation, we
also use $W$ to denote a scalar random variable with the same marginal distribution
$F_W$.}

To  estimate $\theta_0$ from observed data
%\footnote{We  denote by $X_1$, \dots, 
%$X_n$ the rows of $\bX$.
% We often omit
%the arguments $Y$, $\bX$ as this dependency will hold throughout.
%Without loss of generality, we assume that the columns of 
%$\bX$ are normalized so that $\|\bX\,e_i\|_2\approx 1$. (A more
%precise assumption will be formulated below.) }
$(Y,\bX)$ we use an (M)-estimator. Picking  a non-negative
even convex function $\rho : \reals \to \reals_{\ge 0}$, we solve
the optimization problem\footnote{ $X_1$, \dots, 
$X_n$ denote the rows of $\bX$;  while $\theta$ denotes a column vector.
$\htheta$ is chosen arbitrarily if there are multiple  minimizers.
}
\begin{align}
\htheta (Y;\bX)  \equiv \arg\min_{\theta\in\reals^p}  \sum_{i=1}^n
\rho\big(Y_i- X_i \cdot \theta \big)\, ,
\label{eq:Mestimation}
\end{align}
%
%where $\<u,v\>=\sum_{i=1}^mu_iv_i$ is the standard scalar product in
%$\reals^m$, and $\htheta$ is chosen arbitrarily if there are multiple  minimizers.

Of course the prescription is broad enough to encompass traditional
least squares -- $\rho_{LS} (t) = t^2$ -- however, this would not be robust to outliers
\footnote{As can be documented by Frank Hampel's notions of {\sl Influence
  Curve} \cite{hampel1974influence},
which shows that least squares has unbounded influence, and {\sl Breakdown
Point}, which documents  that 
a single bad observation can cause the least squares solution to misbehave 
arbitrarily.}.  Better choices might include least absolute deviations
-- $\rho_{{\rm LAD}}(t) = |t|$ --
and of course the Huber $\rho$ -- $\rho_{{\rm H}}(t; \lambda) = \min( t^2/2 , \lambda |t|-\lambda^2/s)$.\footnote{Other seemingly
good choices, like $\rho(t) = -\log( 1 + t^2)$ are ruled out by lack of convexity.}

\subsection{Fixed $p$, large $n$ Minimax Robustness}

Consider the random design case where $X_i \sim_{\text{iid}} N(0,I_p)$, 
and let  $\psi = \rho'$ denote the \emph{score function}
associated to the $(M)$ estimator of interest.
Let $n \goto  \infty$ with $p$ fixed, and 
consider the per-coordinate asymptotic variance
\[
   V_\infty(\psi,F)  =_{a.s.} \lim_{n \goto \infty}  \frac{n}{p} \cdot \Tr( \Var_F( \hat{\theta} ) ).
\] 
Huber proposed to consider $V_\infty(\psi,F)$ as the payoff function in
a game between the statistician and nature.
The two arguments of $V_\infty$ represent the  two choices  being made here:
the statistician is choosing the estimator, by specifying  $\psi$,
and `nature' is choosing the error distribution, by specifying $F = F_W$.
The statistician pays out
the amount $V_\infty(\psi,F)$ and, planning for all eventualities,
wants to minimize the worst-case payout. The statistician envisions
that  $F$ might contain a fraction $\eps$ of `bad data', and so
assumes that the action space of Nature is 
the class $\cF_\eps$ of all contaminated normal distributions
$F = (1 - \eps) \Phi + \eps H$. Here  $\Phi$ notes the standard normal,
$\eps \in (0,1)$ the contamination fraction,
and $H$ the contaminating distribution. 

For a given choice $\psi$,
the maximal payout that can arise is $\max_{F\in \cF} V_\infty(\psi,F)$.
Huber proposed that the statistician should minimize this quantity
across $\psi$, thus obtaining the minimax asymptotic variance
and the associated minimax score. He found the
least-informative distribution, $F_\eps^*$  \-- the cdf $F$ solving
$\min_{F \in \cF_\eps} I(F)$ with $I$ the Fisher information for location, and 
Huber obtained the formula
\beq \label{eq:classical}
     \min_\psi \max_{F \in \cF_\eps}  V_\infty(\psi,F) = \frac{1}{I(F_\eps^*)}.
\eeq
He also discovered the minimax-optimal score function, now called the Huber score;
it has the form
\[
  \psi_\lambda(x) = \min(\lambda,\max(-\lambda,x)),
\]
for a specific $\kappa=\kappa^*(\eps)$, achieving the minimax.
Numerous textbooks cover this material, including
of course  \cite{huber2009robust}; see also Section \ref{sec-Reminders} below.

\subsection{High-Dimensional Asymptotics}

In his Wald lectures \cite[Page 802]{huber1973robust}
Peter Huber called attention to the fertile regime
beyond the  fixed $p$, large $n$ asymptotic, 
\begin{quotation}
\sl We intend to build an asymptotic theory for $n \goto \infty$; 
but there are several possibilities for the concomitant behavior
of $p$. In particular, with decreasing restrictiveness:
\begin{description}
\item[(a)] $\lim\sup  \, p  <\infty$
\item[(b)] $\lim p^3/n = 0$
\item[(c)] $\lim p^2/n = 0$
\item[(d)] $\lim p/n = 0$
\item[(e)] $\limsup \, p/n < 1$
\item[(f)] $\limsup \,  n - p = \infty$.
\end{description}
\hfill P.J. Huber, {\em Annals of Statistics}, {\bf 1}, 802.
\end{quotation}
Huber also initiated the attack on this hierarchy of new asymptotic settings,
addressing cases (b)-(d).

Though this was 40 years ago, it has taken the profession
a while to catch up\footnote{Peter Bloomfield entered this area already in 1974 \cite{Bloomfield},
and Stephen Portnoy in 1984 \cite{portnoy1984asymptotic}. Soviet-era mathematicians also began studying 
the high-dimensional asymptotic in the late 1960's just when Huber
was also thinking about it; and so Serdobolskii \cite{Serdobolski}
 speaks of the \emph{Kolmogorov asymptotic}, crediting Andrei Kolmogorov
with calculations in the proportional-limit asymptotic already in 1967. Nevertheless,
 Huber's 1972 Wald Lectures were certainly the earliest high-profile venue
 marking out this asymptotic for future research}. In recent years, 
the focus of mathematical statistics research 
has finally gone beyond the fixed $p$, large $n$ asymptotic, to consider regimes (d)-(e),
where  $n$ and $p$ are both large\footnote{A few references here may suffice:
\cite{candes2007dantzig,BickelEtAl,buhlmann2011statistics,karoui2013robust}.}.

In this paper, we consider a precise version of case ($e)$, which
we call the {\em Proportional-Limit } asymptotic $PL(\dof)$;
in this regime $n,p \goto \infty$ and $n/p \goto \dof \in (1,\infty)$.
Thus $\dof$ measures the {\em number of observations per 
parameter to be estimated}.  This parameter
 seems to recur frequently in practitioner thinking:
Huber specifically mentions in his 1972 Wald lectures  the
advice from crystallographers\footnote{Huber's wife Effi Huber-Buser was trained as a 
crystallographer and in the experience of DLD is an insightful scientist,
even knowing quite a lot even about the field of statistics and the statistical profession.} to keep
\footnote{In DLD's first linear models statistics course, based on the classic
Daniel and Wood \cite{daniel1999fitting}, the instructor specifically mentioned $n/p > 10$ as a
desirable ratio. It will be clear from the main results of this paper that 
the prescription to keep $n/p >5$ was very good advice indeed.}  $n/p > 5$.

In this paper the assumption $PL(\dof)$ will further entail a random Gaussian 
design, normalized so for each $n$,
$X_i \sim_{iid} N(0,\frac{1}{n}I_{p\times p})$; and the
regression parameter $\theta_0 = \theta_{0,n} \in \reals^p$ 
will be normalized  so that the per-coordinate size
$p^{-1} \| \theta_{0,n}\|_2  \goto_{a.s.} \tau_0^2$. In this model $ \E\{\bX^{\sT}\bX\} = I_{p \times p}$,
and so under standard Gaussian errors $F=\Phi$,
the per-coordinate Fisher Information is $1$ for every $n$. Because of the finiteness
of the total Fisher Information per coordinate, we are not entitled to
expect highly precise estimation;
hence  it should be no surprise to find that
the MSE  $p^{-1} \| \hat{\theta}_n - \theta_0 \|_2^2  \goto_{a.s.}
\AMSE(\hat{\theta},\theta_0) \neq 0$.
(Here and below $\AMSE$ stands for asymprotic mean square error.)
Consider as performance measure the per-coordinate asymptotic variance:
\[
   V_\dof(\psi,F)  =_{a.s.} \lim_{n \goto \infty}  \frac{1}{p} \cdot \Tr( \Var_F( \hat{\theta}_n ) ).
\]
The notation $V_\dof(\psi,F)$ emphasizes both the dependence
of the asymptotic variance on $\psi$ and $F$ as in the classical case,
but also the dependence on $\dof \in (1,\infty)$. Recent work on (M)-estimates
in $PL(\dof)$ by \cite{karoui2013robust,donoho2013high}
shows that  $V_\dof (\psi,F)  >  V_\infty (\psi,F)$, while
$V_\dof (\psi,F)  \goto V_\infty (\psi,F)$ as $\dof \goto \infty$.

Here we will carry out the Huber program of evaluating  the minimax 
asymptotic variance of the Huber estimate -- this time for $V_\dof(\psi,F)$ for $\dof \in (1,\infty)$,
rather than the classical case $V_\infty$.
The statistician minimizes $V_m$ over
the family $(\psi_\lambda)_{\lambda > 0}$ of all tunings of Huber $(M)$-estimates of regression,
and Nature maximizes over gross-error contaminations  
$F \in \cF_\eps$. 

The classical solution $m=\infty$ plays an important role even
in the $PL(\dof)$ case. Suppose that Huber's least-informative
distribution $F_\eps^*$ obeys
 $I(F_\eps^*)  \cdot \dof > 1$.
In dimensional analysis $I(F_\eps^*)$ is the Fisher information per observation,
while $\dof$ is the number of observations per parameter. Hence this
product is the Fisher information per parameter.
Suppose that this exceeds 1. Then our main result (Corollary \ref{coro-minimax}) shows that
\beq \label{eq:modern}
    \min_\lambda \max_{F \in \cF_\eps}  V_\dof(\psi_\lambda, F) =  \frac{1}{I(F_\eps^*)- 1/\dof} .
\eeq
Of course, the RHS of (\ref{eq:modern})   is strictly bigger than the RHS of 
the classical ($\dof=\infty$) case  (\ref{eq:classical}).
As compared to the classical $\dof=\infty$ case, when $1<\dof < \infty$
 the worst-case asymptotic variance is no longer given by the 
reciprocal of the worst-case Fisher Information. However, the discrepancy  grows
small  as $\dof \goto \infty$.
Hence, new phenomena
emerge in the high-dimensional situation.

\subsection{Variance Breakdown}

Suppose now that the minimal Fisher information per
parameter does not exceed $1$  \-- i.e. that   $I(F_\eps^*) \cdot \dof \leq  1$. Then our main result
additionally states that 
\[
\min_\lambda \max_{F \in \cF_\eps}  V_\dof(\psi_\lambda, F) =  \infty .
\]
In short, the asymptotic variance of the Huber estimator
{\em breaks down} at a critical ratio $\dof = \dof^*(\eps)$ of observations per parameter.
Hampel (1968) defined the breakdown point \--  the minimal
fraction of gross errors that can drive the estimator beyond all bounds.
Later, in connection with non-convex $(M)$ estimators \-- such as Hampel's redescending
(M)-estimator \-- the phenomenon of breakdown of  asymptotic variances arose;
see \cite[Section 5.2]{donoho1983notion}.
For Huber's minimax $(M)$-estimator of classical location, 
no such breakdown occurs: for each $\eps \in (0,1)$,
\[
\min_\lambda \max_{F \in \cF_\eps}  V_\infty(\psi_\lambda, F) < \infty.
\]
 Huber, in personal communication, at one time
considered this non-breakdown of the asymptotic variance 
to be a notable advantage of the Huber estimator
in comparison to some other procedures, 
such as the Hampel `redescending' score function.

Under this paper's $PL(\dof)$ asymptotic, 
variance breakdown of Huber (M)-estimates indeed occurs,
For a fixed ratio $\dof$ of observations per parameter,  the variance
breakdown point is exactly the critical fraction of contamination $\eps$ where the minimal 
Fisher Information per parameter drops to 1 or smaller:
\[
    \eps^* \equiv \eps_\dof^*(\mbox{Minimax Huber-} (M) \mbox{ Estimate}) =  \inf \{ \eps : \dof \cdot I(F_\eps^*) =  1 \}.
\]

\subsection{Illustration}

As a first deliverable of this paper, consider Figure 1,
which displays the minimax asymptotic variance as a function
of the contamination fraction $\eps$ and the degrees of freedom
per parameter estimated $m$.  Below the critical curve -- $1/m =
I(F_\eps^*)$ --
we present contours of the minimax asymptotic variance; in the lower left corner,
the asymptotic variance is nearly $1$, as it would be in the classical $m =\infty$ $\eps=0$ case, 
The minimax asymptotic variance blows up as we approach the dashdot curve.

\begin{figure}[h]
\begin{center}
\includegraphics[height=3in]{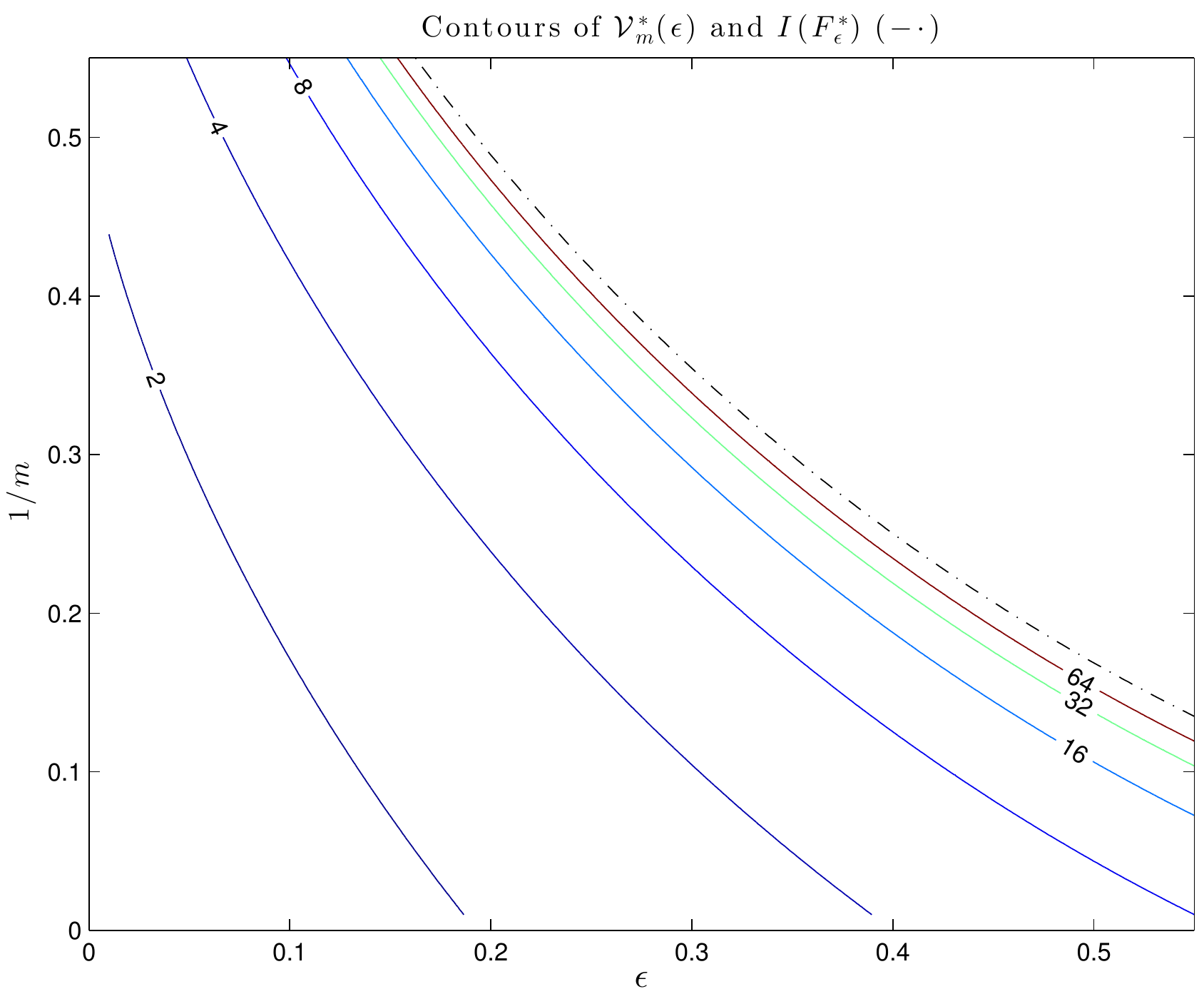}
\caption{Minimax asymptotic variance $\cV^*_m(\eps)$.
Each pair $(\eps,m)$  is represented by the $(x,y)$-point
with $x = \eps$ and $y = 1/m$. 
The resulting parameter space $0 \leq \eps,1/m \leq 1$ is divided into
two phases -- below and above the critical curve
indicated by the dashdot line. Contours of the asymptotic variance
$\cV^*_m(\eps)$ are depicted in the
lower phase; they are undefined in the upper phase, where 
the asymptotic variance cannot be bounded: $\cV_m^*(\eps) = +\infty$.
The boundary separating the 
two phases is indicated by the dashdot curve, at $1/m = I(F_\eps^*)$. }
\label{fig-Fig2}
\end{center}

\end{figure}

\begin{figure}[h]
\begin{center}
\includegraphics[height=3in]{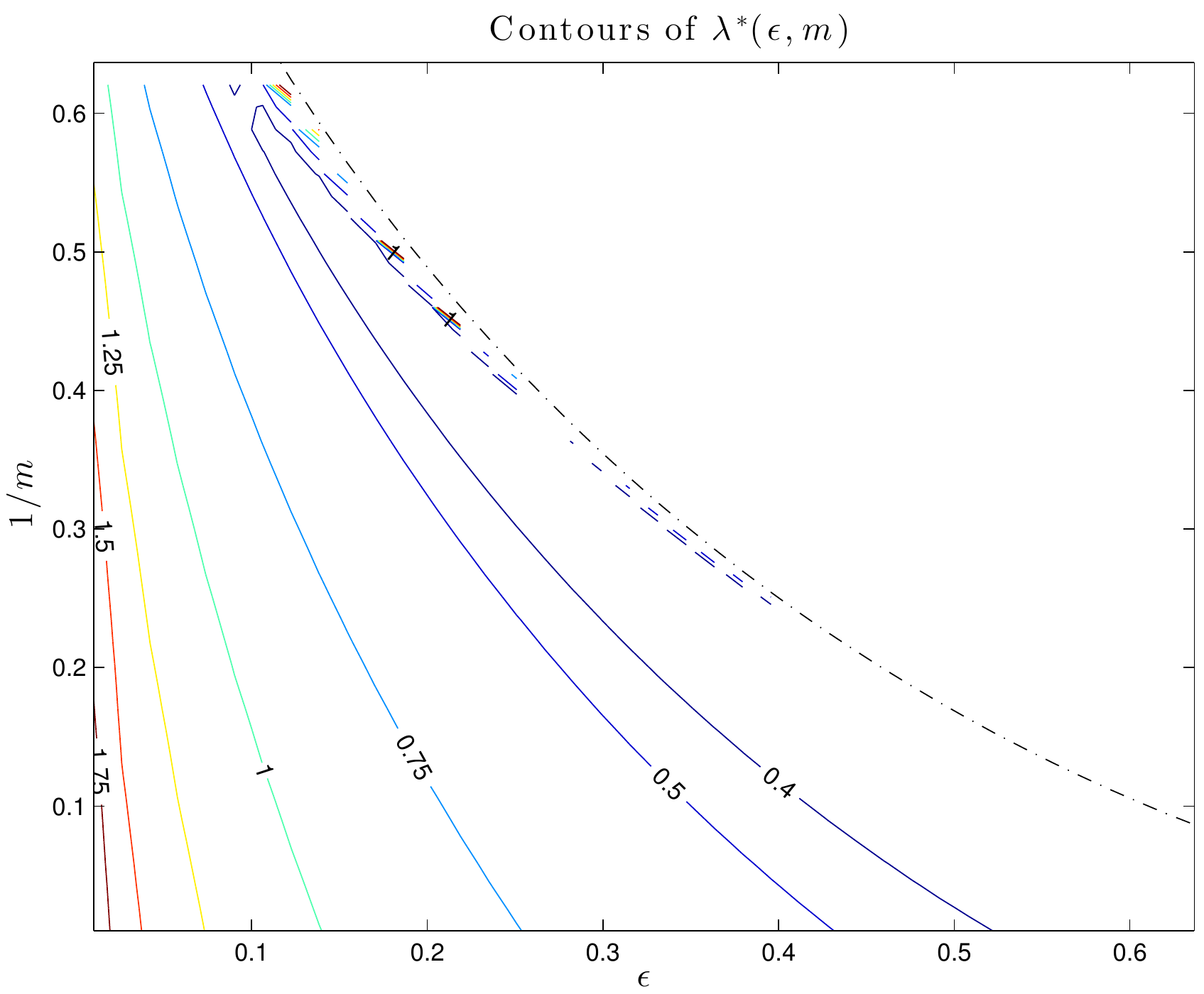}
\caption{Minimax $\lambda^*(\eps;m)$.
Each pair $(\eps,m)$  is represented by the point
$x = \eps$ and $y = 1/m$. 
 Contours of the minimax $\lambda$ parameter
$\lambda^*(\eps;m)$ are depicted in the
region below the dashdot curve
at $1/m = I(F_\eps^*)$. }
\label{fig-Lambda}
\end{center}
\end{figure}

A second deliverable is provided by Figure \ref{fig-Lambda}, which presents contours
of  the minimax tuning parameter
$\lambda^*(\eps,m)$; this selects the Huber $\rho_\lambda$ that
achieves the minimax asymptotic ($V_m$-) variance. Figure \ref{fig-Lambda}
shows that how $\lambda^*$ decays towards zero as $(\eps,m)$ approaches the critical curve.

Table \ref{table:intro} gives some specific numerical values of the minimax asymptotic variance $\cV_m^*(\eps)$.
When $m=2$, it  turns out that the minimax asymptotic variance breaks down at exactly $\eps^* = 0.1924...$,
this is the value of $\eps$ where $I(F_{\eps^*}^*) = 1/2$; the dramatic increase in  variance
as $\eps \uparrow \eps^*$ is plain from the table.

\begin{table}
\begin{center}
\begin{tabular}{r | c c c c c c c}
$\eps$              &  0.05  & 0.10 & 0.15 & 0.175 & 0.1875 & 0.20 & 0.25 \\
\hline
$\cV_2^*(\eps)$ &  3.38 & 5.84 & 13.9 & 35.0 &  136.4 & $\infty$ & $\infty$ \\
\end{tabular}
\caption{Worst-case asymptotic variance of minimax-tuned Huber (M)-estimator, at various levels of contamination;
degrees of freedom per parameter $m=2$.}
\label{table:intro}
\end{center}
\end{table}

We conducted a small Monte-Carlo experiment to illustrate these concepts. With $n = 500$ and $p =250$,
so $m=2$, we considered the linear model with iid Normal predictors $X_{i,j}$, and  contaminated normal errors $W_i$, where
$F_W = G_{\eps,\mu} \equiv (1-\eps) \Phi + \eps H_\mu$, and $H_\mu$ denotes the symmetric Heaviside CDF,
with mass spread equiprobably at $\pm \mu$.  

The reader can see in Table \ref{table:empir} that, 
for small $\eps = 1/20$, even as we make the contamination 
increasingly large, by setting $\mu =100$,
the empirical standard error stays bounded, independently of contamination amplitude $\mu$. However,
as $\eps$ approaches the breakdown point $\eps^*_2 = 0.1924...$, the variance grows considerably as $\mu $ grows large.
%Although the theory shows that with optimal tuning, the Huber estimator still has finite asymptotic 
%variance at $m=2$, $\eps =0.1875$, the empirical behavior at  $n = 500$ is even worse than our 
%asymptotic theory predicts (though this effect would 
%be smaller, in relative terms, at $n=1000$). In short our asymptotic theory 
%seems to understate the extent of blow-up of the empirical variance
%that can arise in finite samples.

%\begin{table}
%\begin{center}
%\begin{tabular}{|r | c c c |}
%\hline
%$\eps$ & $\mu$ & $V_m(\psi_\lambda,G_{\eps,\mu})^{1/2}$ & $\widehat{Var}(\hat{\theta}_n^\lambda)^{1/2}$ \\
%\hline
% 0.05 &    2 &  1.6053 &   1.5883  \\
% 0.05 &    5 &  1.8229 &   1.8662  \\
% 0.05 &   10 &  1.8289 &   1.8801  \\
% 0.05 &   20 &  1.8289 &   1.8594  \\
% 0.05 &  100 &  1.8289 &   1.8436  \\
%\hline
% 0.1875 &    2 &  1.9694 &   1.9900     \\
% 0.1875 &    5 &  3.3404 &   3.5099   \\
% 0.1875 &   10 &  4.4345 &   5.5643  \\
% 0.1875 &   20 &  4.4907 &   8.7302   \\
% 0.1875 &  100 &  4.4907 &  37.8817 \\
% \hline
%\end{tabular}
%\caption{Asymptotic Standard Error and Empirical Standard Error  of minimax-tuned Huber (M)-estimator, at various amplitudes $\mu$ of contamination;
%degrees of freedom per parameter $m=2$. Here the amplitude of the contamination is $\mu$.  For $\eps$ small,
%variability stays controlled as $\mu \goto \infty$, but as $\eps$ approaches the breakdown point $0.1924...$,  variability
%grows very large.}
%\label{table:empir}
%\end{center}
%\end{table}

\begin{table}
\begin{center}
\begin{tabular}{|r | c c |}
\hline
$\eps$ & $\mu$  & $\widehat{Var}(\hat{\theta}_n^\lambda)^{1/2}$ \\
\hline
 0.05 &    2 &   1.5883  \\
 0.05 &    5 &   1.8662  \\
 0.05 &   10 &   1.8801  \\
 0.05 &   20 &   1.8594  \\
 0.05 &  100  &   1.8436  \\
\hline
 0.1875 &    2  &   1.9900     \\
 0.1875 &    5 &   3.5099   \\
 0.1875 &   10  &   5.5643  \\
 0.1875 &   20 &   8.7302   \\
 0.1875 &  100 &  37.8817 \\
 \hline
\end{tabular}
\caption{Empirical Standard Error  of minimax-tuned Huber (M)-estimator, at various amplitudes $\mu$ of contamination;
degrees of freedom per parameter $m=2$. Here the amplitude of the contamination is $\mu$.  These empirical data
reflect this paper's theoretical; conclusion that for  $\eps$ small,
variability stays controlled as $\mu \goto \infty$, but as $\eps$ approaches the breakdown point (here $0.1924...$),  variability
grows very large as $\mu$ increases, even though it will still ultimately stay bounded below the breakdown point.}
\label{table:empir}
\end{center}
\end{table}

\section{Reminders}

\subsection{Classical $(M)$ Estimation and minimax asymptotic variance}
\label{sec-Reminders} 
Huber (1964) supposed we have real scalar observations  $Y_i = \theta_0 + W_i$
where $W_i$ are iid  and symmetrically distributed,
so that $\P ( W > x )= \P (W < -x )$. Hence  $\theta_0 \in \reals$ is the center of symmetry of
the distribution of $Y_i$, and so also  the mean, median, etc.
He introduced the  $(M)$-estimator as a solution $\hat{\theta}$ of
\[
   (M) \qquad  \min_\theta \sum_{i=1}^n \rho(X_i - \theta) ,
\]
where $\rho$ is an even convex function, $\rho(x) = \rho(-x)$, so
the score function $\psi = \rho'$ was monotone nondecreasing.
Under additional regularity conditions, he showed that any solution $\hat{\theta}_n$ obeys
\[
      \sqrt{n} (\hat{\theta}_n - \theta_0) \Longrightarrow_D N(0, V(\psi,F) ), \qquad n \goto \infty,
\]
where the asymptotic variance is given by
\beq \label{eq:HuberAVar}
V(\psi,F) = \frac{ \int \psi^2 \de F}{(\int \psi' \de F)^2}   .
\eeq
For further discussion of regularity conditions, see 
\cite{huber2009robust}.

\begin{figure}[h]
\begin{center}
\includegraphics[height=3in]{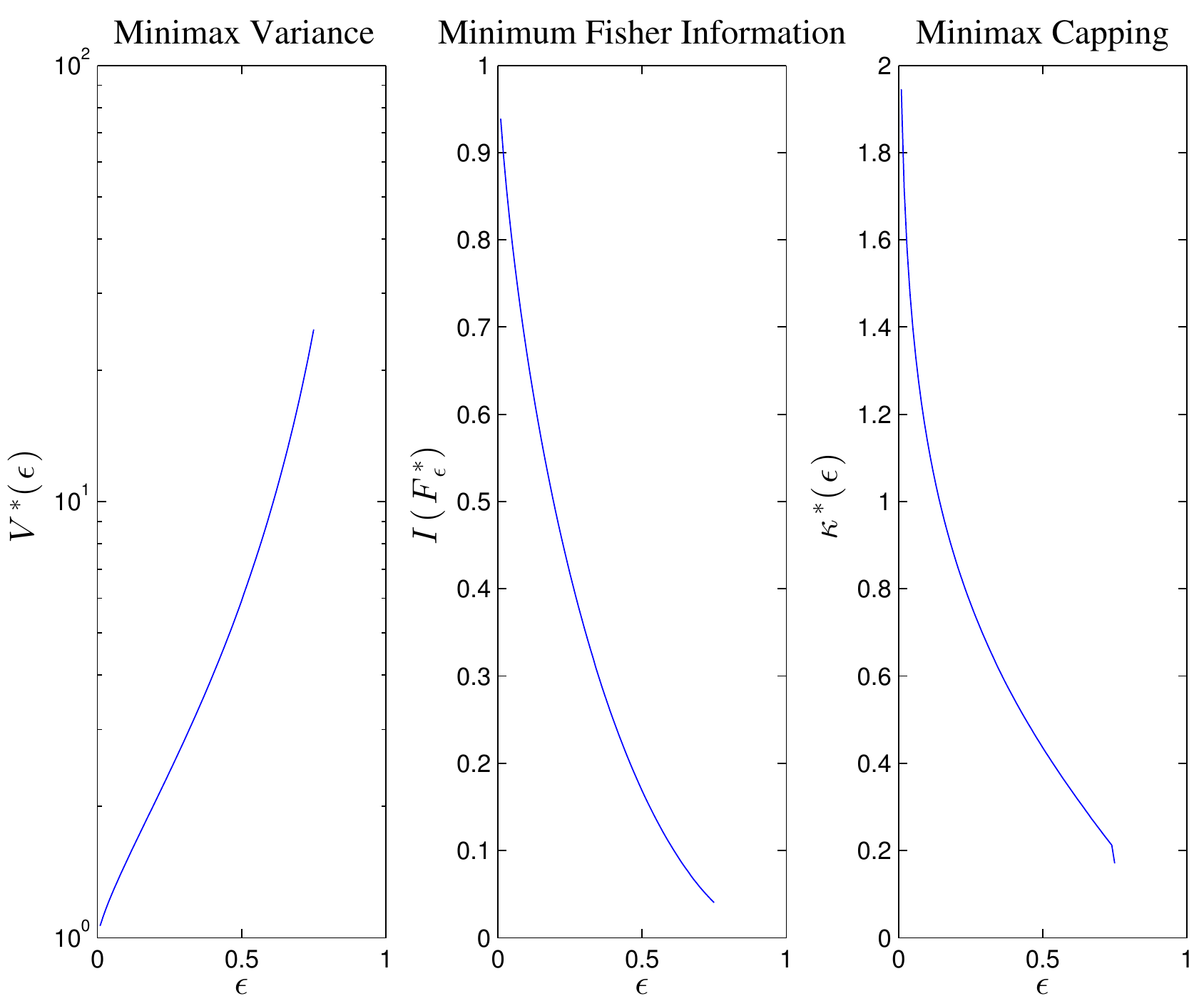}
\caption{Minimax quantities in the \cite{HuberMinimax}
scalar minimax problem, as a function of contamination fraction $\eps$.
From left: $v^*(\eps)$ (semilog plot);
$i^*(\eps)$ and $\kappa^*(\eps)$.
}
\label{fig:HuberMinMaxQuantities}
\end{center}
\end{figure}

Huber considered the situation where the random variable $W_i$ was
distributed roughly as $N(0,1)$, but is subject to gross-errors contamination.
He evaluated
\[
     v^*(\eps)  \equiv     \min_\psi \max_{F \in \cF_\eps}  V(\psi, F) ,
\]
and found the following insightful form. 
Let $I(F) = \int (f'(x))^2/f(x) \de x $ denote the
Fisher information for location;  the
 least informative distribution $F_\eps^*$ 
 minimizes this quantity:
\[ 
     i^*(\eps)  \equiv \min_{F \in \cF_\eps}  I(F) ;
\]
Huber characterized the minimax asymptotic variance 
as the reciprocal of the minimal information:
\[
   v^*(\eps)  = \frac{1}{i^*(\eps)},
\]
and using this was able to write closed formulas
for the optimal shape of $\psi$ -- now called the Huber score function.
In the original paper this was denoted
\[
   \psi_\kappa(x) =  \min(\kappa,\max(x,-\kappa))
\]
with so-called capping parameter $\kappa$, such that errors larger  
in absolute value than $\kappa$
get capped. Huber obtained closed form expressions\footnote{
For example,  $i^*(\eps) = j(\kappa^*(\eps),\eps)$,
where $j(\kappa,\eps) = (1-\eps) \int_{-\kappa}^\kappa x^2 \phi(x) dx + \kappa^2 \cdot \left (\eps + (1-\eps) \cdot 2 \cdot \Phi(-\kappa) \right)$ 
and $\kappa^*(\eps) = \mbox{argmin}_\kappa j(\kappa,\eps).$
}
for the minimax capping parameter $\kappa = \kappa^*(\eps)$,
the least favorable $F = F^*_\eps$, and the minimax asymptotic variance $v^*(\eps) = V(\psi_{\kappa^*(\eps)},F_{\eps}^*)$.
Figure \ref{fig:HuberMinMaxQuantities} displays the behavior of $v^*(\eps)$ and $\kappa^*(\eps)$, as well as $i^*(\eps)$.

\subsection{Regularized Score Functions}

Huber's  $(M)$ estimator of regression uses, for some fixed $\lambda > 0$,
\begin{align}
\rho_{\lambda}(z) = \begin{cases} 
z^2/2 & \mbox{if $|z|\le \lambda$,}\\
\lambda |z| - \lambda^2/2 & \mbox{otherwise.}
\end{cases}\label{eq:HuberDef}
\end{align}
Huber's $\rho$  is quadratic in the middle, has linear tails, and 
is continuous with a continuous derivative. This is straight out of Huber's
theory for the location problem, so no-one should be confused by the switch from $\kappa$ to $\lambda$
to denote the threshold for transition from quadratic to linear; it simply is convenient below to use
$\lambda$ rather than $\kappa$ in the regression case.

For the AMP algorithm discussed below,
we  need the family of {\em regularized} $\rho$-functions,
where for each regularization parameter $r > 0$,
\begin{align}
 \rho(z;r) \equiv\min_{x\in\reals} \Big\{r\rho_{\lambda}(x) +
 \frac{1}{2}(x-z)^2\Big\}\, .
\end{align}
Associated to this  is a 
regularized score function $\Psi(z) = \Psi(z;r)$.
\cite{donoho2013high} writes it in terms 
of Huber's original score $\psi_\lambda$:
\beq \label{eq:relateOrig}
   \dual_{\lambda}( z ; r) =  r \cdot \psi_\lambda \Big( \frac{z}{1+r} \Big) .
\eeq
In particular the shape of each $\Psi$ is similar to $\psi$,
but the slope of the central part is now  $\|\Psi'( \,\cdot\,;r) \|_\infty = \frac{r}{1+r} < 1$.

As explained in \cite{donoho2013high}, although one uses the Huber $\psi$
 as the basis of a high-dimensional regression estimation, the effective score function
 of that $(M)$-estimator belongs to the family $ \dual(\cdot;r)$, for a particular choice of $r$,
 defined below.

\subsection{AMP algorithm}

\newcommand{\res}{R}
The approximate message passing 
(AMP) algorithm we proposed in \cite{donoho2013high}
for the optimization problem (\ref{eq:Mestimation}) 
 is iterative, starting at iteration $0$ with an initial estimate
 $\htheta^0\in\reals^p$. At iteration  $t=0,1,2,\dots$ it 
 applies a simple procedure to update its estimate
 $\htheta^{t}\in\reals^p$, producing $\htheta^{t+1}$.   
 The procedure involves three steps at each iteration.
 \begin{description}
 \item [Adjusted residuals.]
 Using the current estimate $\htheta^t$,
 we  compute the vector of {\em adjusted 
residuals} $\res ^t\in\reals^n$, 
\begin{align}
\res^t & = Y -\bX\htheta^t+\dual(\res^{t-1};r_{t-1})\, ; \label{eq:AMP1}
\end{align}
where to the ordinary residuals $Y - \bX \htheta^t$ we 
here add the extra term\footnote{Here and below, given $f:\reals\to\reals$
and $v=(v_1,\dots,v_m)^{\sT}\in\reals^m$, we define $f(v)\in\reals^m$
by applying $f$ coordinate-wise to $v$, i.e. $f(v) \equiv
(f(v_1),\dots ,f(v_m))^{\sT}$.}
$\dual(\res^{t-1};r_{t-1})$.

\item [Effective Score.]
We choose a  scalar $r_t > 0$,
so that the effective score $\dual(\,\cdot\,;  r_t)$ has empirical average slope $p/n \in (0,1)$.
Setting $\dof=\dof(n) = n/p>1$, we take any solution\footnote{This equation always admits at
  least one solution; cf \cite[Proposition A.1]{donoho2013high}} (for instance the
smallest solution) to
\begin{align}
 \frac{1}{\dof} = \frac{1}{n}\sum_{i=1}^n\dual'(\res^t_i;r)\, .  \label{eq:AMPb}
\end{align}
\item [Scoring.] We apply the effective score function $\dual(\res^t; r_t)$:
\begin{align}
%
%\res^t & = Y -\bX\htheta^t+\dual(\res^{t-1};r_{t-1})\, ,\label{eq:AMP1}\\
\htheta^{t+1} & = \htheta^t+\dof\bX^{\sT}\dual(\res^t;r_t)\, . \label{eq:AMP2}
\end{align}

\end{description}

We emphasize that the above procedure, although presented as an algorithm,
will in fact be used simply a tool in proving results about $(M)$-estimates.

\subsection{State evolution description of AMP}
\label{ssec:StateEvolution}

State Evolution (SE) is a formal procedure  for computing 
the operating characteristics of the AMP iterates $\htheta^t$ and $\res^t$
for arbitrary fixed $t$,  under the $PL(\dof)$ asymptotic  $n,p\to\infty$,
$n/p \goto \dof$.  The ideas have been described at length in \cite{donoho2013high}.
Namely, for the $t$-th iteration of AMP, consider the quantity
\[
  \tau_t^2 \equiv  \lim_{n \goto \infty} \frac{1}{pm}
  \|\hat{\theta}^t-\theta_0\|_2^2 =
  \frac{1}{m}\AMSE(\hat{\theta}^t;\theta_0)\, .
\] 
SE offers a way to calculate $\tau_t$ using $\tau_{t-1}$,
and by extension calculating the limiting AMSE $\dof \lim_{t \goto \infty} \tau_t^2$.

\newcommand{\bu}{{\bf u}}

At the heart of State Evolution 
are  the {\em effective noise level}
$\sigma_t = \sqrt{1 + \tau_t^2}$,
%and the floating threshold 
%$\lambda_t = \kappa \sigma_t$.
which changes iteration by iteration as
the statistical properties of the AMP iterates evolve;
it reflects the combined impact on the estimation of a parameter
of observational noise $W$
with standard deviation $1$ (on the uncontaminated data)
together with estimation noise $\tau$ that `leaks' from the other 
estimated parameters.

Also there is the notion of the  {\em effective
slope}:  the
% for a given $\tau > 0$ and $\dof > 1$ and noise distribution $F_W$,
%let $\sigma = \sqrt{1 + \tau^2} $ and $\lambda = \kappa \sigma$.
%The {\it  effective slope map} is the 
well-defined value  $r = \cR(\tau;\dof,\lambda, F_W)$ giving  the smallest
solution $ r \geq 0 $ to
\[
\frac{1}{\dof} =   \E\Big\{\dual'_{\lambda} (W+\tau\, Z;r)\Big\},
\]
where $W\sim F_W$, and, independently, $Z\sim N(0,1)$.
Informally, $\cR$ measures the value of the  regularization parameter $r$
that  satisfies the population analog of the AMP empirical average slope condition (\ref{eq:AMPb}).

Similarly, define the \emph{variance map}
\[
   \cA(\tau^2 , r;  \lambda, F_W)=   \E\Big\{\dual^2_{\lambda}(W+\tau\, Z;r)\Big\}\, ,
\]
$\cA$ measures the variance of the
resulting effective score. Evidently, for $r > 0$, $0 \leq \cA(\tau^2, r)  \leq (\Var(W)+\tau^2)$.

In the last two displays, the reader can see that extra Gaussian noise of
variance $\tau^2$ is being added to the underlying noise $W$. 
%We emphasize that
%$\lambda$ is in both cases floating with $\sigma$ according to $\lambda = \kappa \sigma$.

\newcommand{\tcV}{\tilde{\cal V}}
\newcommand{\bcV}{\overline{\cal V}}
\begin{definition}
{\em State Evolution} is an iterative process for computing the sequence of scalars
$\{\tau^2_t\}_{t\ge  0}$,  starting from an initial condition 
$\tau_0^2\in \reals_{\ge 0}$ following the recursion
\begin{align}
\tau_{t+1}^2 =  \dof \cdot \cA(\tau_t^2 , \cR(\tau_t) ) =  \dof \cdot \cA(\tau_t ^2, \cR(\tau_t; \dof,\lambda,F_W) ; \kappa, F_W) .\label{eq:StateEvolution}
\end{align}
\end{definition}

Defining  $\cT(\tau^2)=  \dof \cdot \cA (\tau^2 , \cR(\tau) )$, we see that the evolution of $\tau_t^2$
follows the iterations of the map $\cT$. In particular, we make these
observations:
\bitem
 \item $\cT(0)  > 0 $,
 \item $\cT(\tau^2)$ is a continuous, nondecreasing function of $\tau$.
% \item $\tcV(\tau^2)$ is a concave function of $\tau^2$.
 \item $\cT(\tau^2) < c\cdot \tau^2$ for some $c\in (0,1)$ and all sufficiently large  $\tau$.
\eitem

As a consequence of Theorem \ref{thm:SECorrectHuber} below,
$\cT$ has a unique 
 fixed point $\tau_\infty^2$, i.e. 
 \[
     \cT(\tau_\infty^2) = \tau_\infty^2 .
 \]
If follows from the above properties that  this fixed point is stable
and attracts $(\tau_t^2)$ from any starting value. Explicitly,
   for each initial value $   \tau_0 \in (0,\infty) $, the sequence defined for $t =1,2,\dots$ by
  $\tau_t^2 = \cT(\tau_{t-1}^2)$ converges to the above fixed point:
 \[  
       \tau_t^2 \goto \tau_\infty^2, \qquad  \mbox{ as  } t\goto \infty  .
 \]
 \subsection{Correctness of State Evolution}
 
 The paper \cite{donoho2013high} 
 considers $(M)$ estimates with strongly convex $\psi$-functions --
 this excludes the Huber estimator for technical reasons.
 In that paper,  \cite[Theorem 3.1]{donoho2013high} shows that State Evolution correctly computes the
 operating characteristics of the AMP algorithm.
 In particular, the AMP algorithm has $\dof \cdot \tau_\infty^2$ for its $t \goto \infty$
limiting AMSE in estimating $\theta_0$.

Within the strongly convex setting, \cite[Theorem 4.1]{donoho2013high} shows that the AMP algorithm
converges in mean square to the (M)-estimator, which is therefore also 
described by the fixed point of State Evolution. 

Define the asymptotic variance of  the $(M)$-estimator $\htheta$
by
\[
  \AVar(\htheta) = \lim_{n,p \to \infty} {\rm Ave}_{i\in[p]}\Var(\htheta_i),
\]
where ${\rm Ave}_{i\in[p]}$ denotes the average across indices $i$.
\cite[Corollary 4.2]{donoho2013high} shows that
the asymptotic variance of $\htheta$ obeys
\begin{align}
\AVar(\htheta_i)= \dof  \tau^2_\infty . 
\end{align}

It follows that State Evolution describes not only the operating characteristics
of the large $t$-limit of the AMP algorithm, but {\em any} algorithm for obtaining the (M)-estimate
in the $PL(\dof)$ asymptotic.
So the fixed point of the one-dimensional dynamical system $\tau^2 \mapsto \cT(\tau^2)$
is fundamental.

All these results extend to the Huber estimator itself. The companion paper
\cite{andreaHuberRigor} proves the following
extension of the results in \cite{donoho2013high}.
\begin{thm} \label{thm:SECorrectHuber}\cite{andreaHuberRigor} 
Suppose that $(\tau_\infty,r_\infty)$ solve the two equations
\begin{eqnarray*}
\frac{1}{\dof} &=&   \E\Big\{\dual'_{\lambda} (W+\tau\, Z;r)\Big\},\\
\tau^2 &=&   m\E\Big\{\dual^2_{\lambda}(W+\tau\, Z;r)\Big\}\, .
\end{eqnarray*}
Then under the $PL(m)$-limit, 
 the Huber $(M)$-estimator $\hat{\theta} = (\hat{\theta}_i)_i$ obeys:
\[
     \AVar(\hat{\theta}) = m \cdot  \tau_\infty^2 .  
\]
In particular, this implies that such fixed point is unique.
\end{thm}

We note that, at the fixed point $(\tau_\infty,r_\infty)$, we have
\[
 \AVar(\hat{\theta})  =  \frac{\E\Big\{\dual^2_{\lambda}(W+\tau\,
   Z;r)\Big\}}{\E\Big\{\dual'_{\lambda} (W+\tau\, Z;r)\Big\}^2}\, .
\]
the expression on the RHS can be written in terms of Huber's asymptotic variance formula $V$ (\ref{eq:HuberAVar}): it is
 $V(\dual_{\lambda}(\cdot ;r), F_W \star N(0,\tau_\infty^2)) $.
In other words, the classical Huber asymptotic variance formula continues to hold in an extended sense; however, it is evaluated at
the effective score function $\dual_{\lambda}(\cdot ;r)$ with respect to the effective error distribution
$ F_W \star N(0,\tau_\infty^2)$; see \cite{karoui2013robust} for another approach to this formula.

\newcommand{\ttau}{\bar{\tau}}
\newcommand{\tsigma}{\bar{\sigma}}
\newcommand{\tlambda}{\bar{\lambda}}
\newcommand{\btau}{\bar{\tau}}
\newcommand{\bbee}{\tilde{b}}
\newcommand{\bsigma}{\tilde{\sigma}}
\newcommand{\blambda}{\bar{\lambda}}

\section{Least-Favorable State Evolution}

In this section we 
develop an upper bound on the behavior of State Evolution.
We first introduce a variant of SE, in which $\lambda$ evolves rather than staying fixed.
This variant can be conveniently analyzed. In a later section,
we tie the results obtained for this evolution to the original state evolution.

\subsection{Floating-Threshold State Evolution}
\label{ssec:FTSE}

Recall the  notion of {\em effective noise level}
$\sigma_t = \sqrt{1 + \tau_t^2}$ in state evolution,
and consider a variant of SE
where the threshold parameter $\lambda_t$ 
`floats' proportionally to the noise level $\sigma_t$,
as follows $\lambda_t = \kappa \cdot \sigma_t$.
Here $\kappa$ may be viewed as the capping parameter for data which are presumed to be standardized,
and so  the floating  $\lambda_t$  is actually invariant across iteration -- when expressed in
multiples $\kappa$ of the effective noise level.

In an abuse of notation, define $\cA$ with a $\kappa$ (rather than $\lambda$)
as argument to be the variance map, based on floating $\lambda$:
\[
   \cA(\tau^2 , r;  \kappa, F_W)=   \E\Big\{\dual^2_{\kappa \cdot \sigma}(W+\tau\, Z;r)\Big\}\, .
\]
Compounding the abuse, define  $r = \cR(\tau;\dof,\kappa, F_W)$ analogously, so
that
\[
\frac{1}{\dof} =   \E\Big\{\dual'_{\kappa \cdot \sigma} (W+\tau\, Z;r)\Big\}.
\]
Similarly, we define $\cT(\tau^2; m, \kappa, F_W) =  \dof \cdot \cA (\tau^2 , \cR(\tau; m, \kappa, F_W) )$, without any warning
to the reader that the same symbols are being used as in the earlier state evolution with fixed $\lambda$ while
here and below the appearance of $\kappa$
in the argument always refers to the floating $\lambda$ evolution. 
For example, we might write $\tau_\infty^2(m,\kappa,F)$ for the fixed point of a floating-$\lambda$ evolution
and $\tau_\infty^2(m,\lambda,F)$ for the (in general different) fixed point of a fixed-$\lambda$ evolution.
As a first justification for this, note that
the fixed points of the two different dynamical systems (fixed- and floating- $\lambda$ dynamical systems) are in
one-one correspondence, via
\[
     \lambda = \kappa \cdot \sqrt{1 + \tau_\infty^2},
\]
i.e.
\bitem
\item  The fixed-$\lambda$  fixed point  $\tau_\infty^2(\lambda)$ is identical to 
the floating-$\lambda$ fixed point $\tau_\infty^2(\kappa_\infty(\lambda))  $,
under the floating-$\lambda$ parameter 
$  \kappa_\infty (\lambda) = \lambda / \sqrt{1 + \tau_\infty^2(\lambda) }$;  while 
\item The floating-$\lambda$  fixed-point
$\tau_\infty^2(\kappa)$ is identical to the  fixed-$\lambda$  
fixed point
$\tau_\infty^2(\lambda_\infty(\kappa)) $    
at parameter  $\lambda_\infty =  \kappa \cdot \sqrt{1 + \tau_\infty^2(\kappa) }$ .
\eitem
Setting $\lambda = \lambda_\infty(m,\kappa,F_W)$ and $\kappa = \kappa_\infty(m,\lambda, F_W)$
establishes the correspondence. 
%Each of the evolutions has a fixed point at the corresponding situation.
Hence characterizing the fixed points of the floating $\lambda$ scheme will also
characterize those of the fixed lambda scheme; see also Definition \ref{def:calibrate} et seq. below.

\subsection{Least-Favorable SE}
\label{ssec:LFSE}
%We are now in a position to explain the relationship
%of WCSE introduced in  section \ref{ssec:WCSE}
%as the least favorable variant of the SE procedure
%of the previous section.

%For each $\sigma \geq 1$,
% let $\cF_{\eps,\sigma}$ be the collection 
%of gross-error contaminations of 
%a normal distribution with standard deviation $\sigma$:
%\[
% \cF_{\eps,\sigma} = \{F:  F = (1-\eps) \Phi_\sigma + \eps H \},
%\]
%so that of course $\cF_{\eps,1}$ is the same as what we earlier called $\cF_\eps$.

Let $H_{\infty}$ denote the improper distribution 
with its probability mass placed evenly on $\{ \pm \infty \}$; with this notation,  set
 $\bar{F}_{\eps} = (1-\eps) \Phi + \eps H_\infty$.  We now describe an extremal form of floating-threshold 
 state evolution.

\begin{definition} \label{def-LFSE}
{\em Least Favorable State Evolution} (LFSE) is an iterative process for computing a sequence of scalars
$\{\ttau^2_t\}_{t\ge  0}$,  starting from an initial condition 
$\ttau_0^2\in \reals_{\ge 0}$. An instance of LFSE is
determined by $\tau_0^2$ together with fixed positive scalars $\dof$, $\kappa$ and $\eps$.
%Consider now a seemingly different notion of state evolution, based on the
%family $\bar{X}_{\eps,\sigma}$. 
%Let $\ttau_0 = \btau_0$ and define
%the sequence of $\ttau_t$ for $t=1,2,\dots$.

At  the $t$-th iteration,
one needs the $(t-1)$'th result $\ttau_{t-1}$ and
sets $\tsigma_{t-1} = \sqrt{1 + \ttau_{t-1}^2}$ , 
\begin{eqnarray*}
\tarr_t &=&  \cR(\ttau_{t-1};\dof,\kappa,\bar{F}_{\eps}) \\
\ttau_t^2 &=& \dof \cdot \cA(\ttau_{t-1} ^2, \tarr_t ; \kappa, \bar{F}_{\eps})  
\end{eqnarray*}
The procedure is then repeated at the next iteration $t+1$, and so on
\end{definition}
 
Letting $\Phi_\sigma$ denote the CDF for $N(0,\sigma^2)$,
 set $\bar{F}_{\eps,\sigma} = (1-\eps) \Phi_\sigma + \eps H_\infty$,
and define an
improper random variable  
 $\bar{X}_{\eps,\sigma} \sim \bar{F}_{\eps,\sigma}$, taking infinite values
with positive probability. 
Setting $\sigma = (1 + \tau^2)^{1/2}$,
we have $ \bar{F}_{\eps,\sigma}  = \bar{F}_{\eps}  \star  \Phi_{\tau}$.
Definition \ref{def-LFSE}, written in terms of  the improper random variable
$\bar{X}_{\eps,\tsigma_{t-1}}$, and the floating threshold $\tlambda_t= \kappa \cdot \tsigma_{t-1}$, gives:
\[
    \frac{1}{\dof} = \E \Psi'_{\tlambda_t} (
    \bar{X}_{\eps,\tsigma_{t-1}}, \tarr_t) \, ,
\]
and  
\[
    \ttau^2_{t}  = \dof \cdot \E \Psi_{\tlambda_t}^2 (
    \bar{X}_{\eps,\tsigma_{t-1}}, \tarr_t)\, .
\] 

Although $\bar{X}_{\eps,\tsigma_{t-1}}$ is an improper random variable,
these expectations are well defined\footnote{given the boundedness and differentiability of the  underlying Huber $\psi$}.
We refer to instances where state
evolution is applied to proper distributions  in $\cF_{\eps}$
as {\em proper state evolutions}.

\begin{lemma}\label{lem:dominance}{\bf  (LFSE Dominates.)}
Consider a given instance $(\dof,\tau_0,\kappa,F)$ of floating-threshold state evolution 
where $F \in \cF_\eps$.
The LFSE instance $(\dof,\tau_0,\kappa,\eps)$ 
dominates this proper state evolution, namely:
 with $\tarr_t$ the sequence of LFSE regularizing parameters 
 and $r_t$ the sequence of proper SE regularizing parameters,
\[
      \tarr_{t}\geq r_t\, , \qquad t=1,2,\dots ,
\]
while for $\ttau_t^2$ the MSE under LFSE   
and $\tau_t^2$ under proper SE, respectively, we have:
\[
    \ttau_{t}^2 \geq {\tau}_t^2 \, , \qquad t=0,1,2,\dots
\]
\end{lemma}

\begin{figure}[h]
\begin{center}
\includegraphics[height=3in]{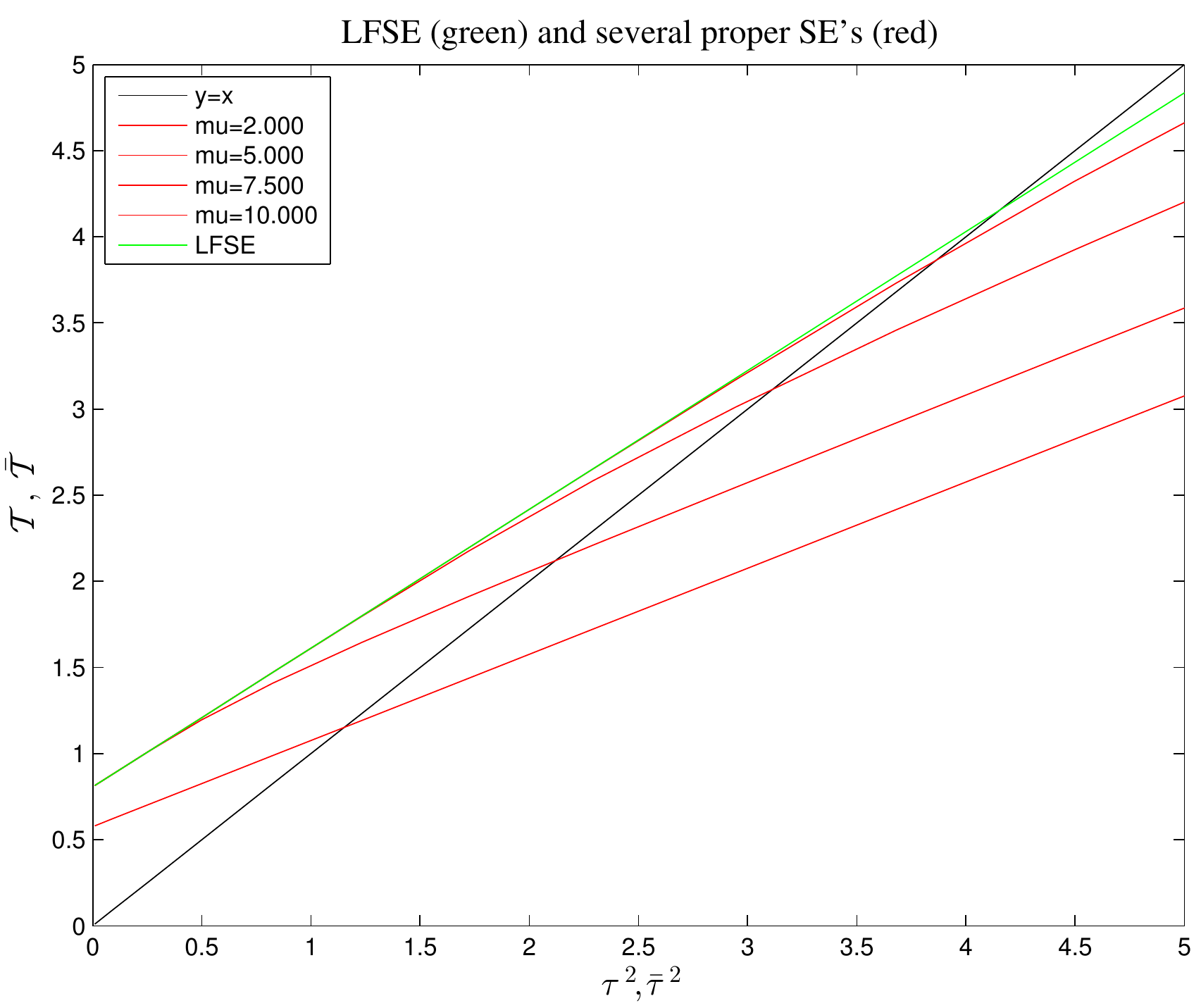}
\caption{MSE maps of proper state evolutions  and of LFSE.
Here $\eps = 0.05$, $m = 5$, and $\mu = 2,5,7.5,10$. The variance map of LFSE
is the green straight line, which lies above the variance maps of all the
proper SE's as depicted by red curves.  Correspondingly, its
fixed point is also higher.}
\label{fig-BoundedPhaseSE}
\end{center}
\end{figure}

Figure \ref{fig-BoundedPhaseSE} illustrates the dominance of LFSE; it shows that the
corresponding dynamical maps obey $\bar{\cT} \geq \cT$.
 
The proof  \-- given in the appendix \--   will depend on the following sequence of observations:
\begin{lemma} \label{lem:monotonicity}{\bf Monotonicity in $x$, $r$, and $\lambda$.}
Let $\Psi( x,r)$ denote the regularized score function based on
Huber's $\psi_\lambda$.
(With $\lambda$ fixed unless stated otherwise.)
\begin{enumerate}
 \item For each fixed $r \in \reals_+ $, $\Psi_{\lambda}( x ,r)$ is a monotone increasing function  of $|x|$;
 \item For each fixed $r \in \reals_+ $, $ \Psi'_{\lambda}( x ,r)$ is a monotone nonincreasing function  of $|x|$;
 \item  For each fixed $x \in \reals$; $r \mapsto |\Psi_{\lambda} (x,r)|$  is monotone nondecreasing in $r$; and
 \item For each fixed $x \in \reals$, $\lambda \mapsto | \Psi_{\lambda} (x,r)|$  is monotone nondecreasing in $\lambda$.
\end{enumerate}
\end{lemma}

It will also need the following invariances,
which  are very special to the extremal improper RV's $\bar{X}_{\eps,\tsigma}$ and $\bar{X}_{\eps,{\sigma}}$ together with
the fact that the proper SE and LFSE use exactly the same $\kappa$ in forming their respective floating $\lambda$'s.

\begin{lemma} \label{lem:invariances}
For $r > 0$, and $t \geq 1$, let $0 < \sigma_{t-1}  < \tsigma_{t-1}$
and $\lambda_t = \kappa \sigma_{t-1}$, and $\tlambda_t = \kappa \cdot \tsigma_{t-1}$.
\beq \label{eq:invar}
   \E \Psi'_{{\lambda}_t} (\bar{X}_{\eps,{\sigma}_{t-1}}, r)  =  \E \Psi'_{\tlambda_t} (\bar{X}_{\eps,\tsigma_{t-1}}, r); 
\eeq
\beq \label{eq:quad-invar}
   \E \Psi_{\lambda_t}^2(X_{\eps,\sigma_{t-1}},r) = \left( \frac{\sigma_{t-1}}{\tsigma_{t-1}} \right)^2 \cdot  \E \Psi_{\tlambda_t}^2(X_{\eps,\tsigma_{t-1}},r).
\eeq

\end{lemma}

\newcommand{\Bbar}{\bar{B}}
\subsection{The  envelope functionals  $\bar{A}$ and $\bar{B}$}

To make LFSE
more  transparent,
we introduce some helpful notation.
In this subsection, we are again in Huber's original
location setting. The evaluation of $v^*(\eps)$
is made significantly easier by helpful notation. Suppose that
$F$ is a sub distribution, i.e. a CDF on the extended reals, and put
\begin{eqnarray}
A(\psi_\kappa,F)  &=& \int_{-\infty}^\infty \psi_\kappa^2(w) \de F(w) = \E_F  \psi_\kappa^2(W) \\
B(\psi_\kappa,F) &=&  \int_{-\infty}^\infty \psi'_\kappa(w) \de F(w) = \E_F  \psi_\kappa'(W) 
\end{eqnarray}
where $W \sim F$. Calculating explicitly for the Huber score function,
we can equally well write
\[
  A(\psi_\kappa,F) = \int_{-\kappa}^{\kappa} w^2\de F(w)  + \kappa^2 \cdot \P_F \{ |W| \geq \kappa \}.
\]
and
\[
  B(\psi_\kappa,F) = \P_F \{ |W| \leq \kappa \}.
\]
%
% Not needed, and incorrect:
%
%%We will use  two (obvious) invariances: if $\sigma^2 = (1 + \tau^2)$ then
%\begin{eqnarray}
%A(\psi_{\kappa \cdot \sigma}, F \star N(0,\tau^2))  &=&  A(\psi_\kappa, F)  \cdot \sigma^2  \\
%B(\psi_{\kappa \cdot \sigma} ,F \star N(0,\tau^2)) &=&   B(\psi_\kappa, F)   .
%\end{eqnarray}

Now define the envelope functions $\bar{A}$ and $\Bbar$, so that
\begin{eqnarray}
\bar{A}(\kappa,\eps)  &=&  \sup \{ A(\psi_\kappa, F) : F \in \cF_\eps  \}  \\
\Bbar(\kappa,\eps) &=&   \inf \{ B(\psi_\kappa, F) : F \in \cF_\eps \}  .
\end{eqnarray}
More explicitly, with $\Phi$ denoting the standard normal CDF, and $Z \sim N(0,1)$,
\begin{eqnarray}
\bar{A}(\kappa,\eps) &=& (1-\eps) A(\psi_\kappa,\Phi) + \eps \kappa^2 \\
\Bbar(\kappa,\eps) &=&   \P_F \{ |Z| \leq \kappa \} = 2\Phi(-|\kappa|)-1 .
\end{eqnarray}
Defining $V(\psi_\kappa,F) = A(\psi_\kappa,F)/B^2(\psi_\kappa,F)$, and correspondingly,
\[
\bar{V}(\kappa,\eps)  =   \sup \{ V(\psi_\kappa, F) : F \in \cF_\eps  \} .
\]
It follows from Huber(1964) that
\[
\bar{V}(\kappa,\eps)  =   \frac{\bar{A}(\kappa,\eps)}{\Bbar^2(\kappa,\eps)} ,
\]
(the inequality LHS $\leq$ RHS is obvious) and also that 
\[
v^*(\eps) = \inf_\kappa \bar{V}(\kappa,\eps) ;
\]
(the inequality LHS $\leq$ RHS again being immediate).

\subsection{Explicit Solution of Least Favorable State Evolution}
\label{ssec:ALT-LFSE}

We now put Huber's notation from the previous subsection
to work, giving explicit formulas for LFSE. 
\begin{lemma}
For a given tuple $(\dof, \eps, \kappa)$  obeying $(1-\eps) >  1/m$,
there is a unique positive solution   $\bbarr(\dof, \eps, \kappa)$  to
\begin{equation} \label{eq:rbardef}
\left( \frac{\bbarr}{1+\bbarr} \right) \cdot \Bbar(\kappa \cdot (1 + \bbarr) ,\eps)  = \frac{1}{\dof} .
\end{equation}
\end{lemma}

Using this notation, we give an explicit characterization 
of LFSE. Let $\bkappa = \bkappa(m,\eps,\kappa) = \kappa \cdot (1 + \bbarr)$ as in
the first argument of $\bar{B}$ in (\ref{eq:rbardef}).

\begin{lemma} \label{lem:LFSEFP}
LFSE with parameters $(m,\eps,\kappa)$ satisfies, with $\bkappa = \bkappa(m,\eps,\kappa)$
\[
  \bcT(\tau^2; \dof,\eps,\kappa) =  (1 + \tau^2) \cdot \bar{V}( \bkappa , \eps ) / m , 
\]
and, if $\bar{V}( \bkappa , \eps ) < \dof$, LFSE has the unique stable fixed point
\[
  \bar{\tau}^2_\infty(m,\eps,\kappa) = \frac{\bar{V}( \bkappa , \eps ) /m}{1-\bar{V}( \bkappa , \eps )/m } .
\]
\end{lemma}

To prove this, consider a seemingly different evolution,
which we call {\em double-bar evolution}:
with $\bbarr$ as introduced above, define
\beq \label{eq-dubbar}
    \bbcA(\tau^2; m,\eps,\kappa ) =   (1+ \tau^2)  \cdot \bar{V}( \bkappa , \eps ) /m^2 % \left( \frac{\bbarr}{1+\bbarr} \right)^2 \cdot  \bar{A}( \bkappa , \eps ) 
\eeq
and
\[
\bbcT(\tau^2) = \dof \cdot {\bbcA}( \tau^2).
\]

With $\dof$,$\eps$, $\kappa$ and thus $\bbarr$ and $\bkappa$  fixed,
define a sequence $\bbtau_t^2$ for $t = 0,1,2,\dots$.
%and an associated sequence
%\[
%   \lambda_t = \kappa \sqrt{1 + \btau_t^2} , \qquad t = 0,1,2,\dots .
%\]
At iteration $t=0$, we  pick a starting value $\bbtau_0 \geq 0$,
we then proceed inductively,
%and let $\lambda_0 = \kappa \sqrt{1 + \btau_0^2}$.  
setting  all later iterates by :
\[
     \bbtau^2_{t} = \bbcT(\bbtau_{t-1}^2),  \qquad t =1,2,\dots.
\]

%Double-bar evolution has its regularization parameter 
%$\bbarr_t = \bbarr$, fixed independently 
%of iteration counter $t$. Consequently, $\bkappa$ is likewise fixed, and  the coefficient
%of $(1 + \tau^2)$ on the RHS of (\ref{eq-dubbar} is fixed throughout the evolution.
%By (\ref{eq:rbardef}) and the definition of $\bar{V}$, the coefficient can be rewritten
%\[
%\bar{V}( \bkappa , \eps ) /m =  \left( \frac{\bbarr}{1+\bbarr} \right)^2 \cdot  \bar{A}( \bkappa , \eps ) 
%\]
Now (\ref{eq-dubbar})  
sets up the dynamical system $\bbtau^2 \mapsto \bbcT(\bbtau^2)$
as an {\em affine} dynamical system (in the variable $\bbtau^2$). Its fixed point (if it exists at all) 
must obey
\[
 {\bbtau}_\infty^2 = (1 +  {\bbtau}_\infty^2) \cdot \bar{V}( \bkappa , \eps ) /m.
\]
So  double-bar evolution has the following explicit solution:

\begin{lemma}
Consider the  double-bar evolution introduced in this section,
with parameters $(m,\eps,\kappa)$. If $\bar{V}( \bkappa , \eps ) < \dof$, it has the unique stable fixed point
\[
  {\bbtau}^2_\infty(m,\eps,\kappa) = \frac{\bar{V}( \bkappa , \eps ) /m}{1-\bar{V}( \bkappa , \eps )/m } .
\]
Otherwise there is no fixed point, and successive iterates run off to infinity.
\end{lemma}

In fact, double-bar evolution is really just LFSE, in disguise. Results of the next subsection
will prove:

\begin{lemma} \label{lem:identSE}
With $(\tarr_t)$ and $(\ttau_t)$ defined by the  procedure of Section \ref{ssec:LFSE},
and  $(\bbarr_t)$  $(\bbtau_t)$ defined by the procedure of this section, 
each initialized identically -- $\ttau_0 = \bbtau_0$ -- we have
\[
    \bbarr_{t} = \tarr , \qquad t=0,1,2,\dots,
\]
\[
    \ttau_{t} = \bbtau_t , \qquad t=0,1,2,\dots,
\]
and
\[
  \bbcT (\cdot) =  \bcT (\cdot).
\]
\end{lemma}

Lemma \ref{lem:LFSEFP} then follows from the last two lemmas.
In turn, Lemma \ref{lem:identSE} follows immediately from the following:
\begin{lemma} \label{lem:extremalSE}
\begin{equation} \label{eq:max-r}
 \sup_{F \in \cF_{\eps} } \cR(\tau; m,\kappa,F  ) = \bbarr(m,\eps,\kappa),
\end{equation}
\begin{equation} \label{eq:max-A}
 \sup_{F \in \cF_{\eps} } \cA( \tau^2, \cR(\tau;m,\kappa,F ) ; m,\kappa, F)   = \bbcA(\tau^2;m,\eps,\kappa).Œ
\end{equation}
\end{lemma}
This shows that  the affine evolution (\ref{eq-dubbar}) 
indeed implements LFSE, and proves Lemma \ref{lem:identSE}.

The proof of Lemma \ref{lem:extremalSE} is given in the Appendix; it depends
on terminology and results of the next subsection.

\subsection{Bounds for $\cA$}

The quantity $\cA$ 
occurring in LFSE is defined using
moments of $\E \Psi^2$; however, Section \ref{ssec:ALT-LFSE}
defines $\bbcA$ in terms of  $\bar{A}$,
which uses moments of $\psi^2$.
To explain the connection -- and prove Lemma \ref{lem:identSE} --
we need to relate the two kinds of moments.

Indeed,  (\ref{eq:relateOrig}) says that 
\[
 \dual( z ; r) =  \frac{r}{1+r} \psi_{\lambda (1+r) }\Big( {z}\Big)\, ,
 \]
 and we also have, for any random variable $X$,
 \[
    \E \psi_\lambda^2(c X ) = c^2 \E \psi_{\lambda/c} ^2 (X),
 \]
 while
 \[
    \E \psi'_{ \lambda}(c X  ) =  \E \psi'_{ \lambda/c}(X). 
 \]
Furthermore, supposing that $W$ has distribution $(1-\eps) \Phi + \eps H$
and that $Z \sim \Phi$ while $U \sim H$, then
\[
   \E |\Psi(W + \tau Z ; r)|^2 = (1-\eps) \E |\Psi(\sqrt{1+\tau^2} \cdot  Z ; r)|^2  + \eps \cdot 
                                                                \E |\Psi(U + \tau Z ; r)|^2.
\]
Now introducing $a \equiv \sqrt{1+\tau^2}/(1+r)$  where $r$ is some fixed positive scalar
kept the same in all the coming displays,
\[
\E |\Psi(\sqrt{1+\tau^2} \cdot  Z ; r)|^2 =  \left( ar \right)^2 \cdot  \E \psi^2_{\lambda/a}(Z) =  \left(  ar \right)^2 A(\lambda/a,\Phi).
\]
and
\[
\E |\Psi(U + \tau \cdot  Z ; r)|^2 =  \left( \frac{r}{1+r} \right)^2 \cdot  \E \psi^2_{\lambda \cdot (1+r) }(U+\tau Z) 
 =  \left( \frac{r}{1+r} \right)^2 A(\lambda \cdot (1+r) , H \star \Phi_\tau).
\]
Similarly,
\[
   \E \Psi'(W + \tau Z ; r) = (1-\eps) \E \Psi'(\sqrt{1+\tau^2} \cdot  Z ; r)  + \eps \cdot 
                                                                \E \Psi'(U + \tau Z ; r);
\]
and
\[
   \E \Psi'(\sqrt{1+\tau^2} \cdot  Z ; r)  = \left( \frac{r}{1+r} \right) \cdot \E \psi'_{\lambda (1+r)} ( \sqrt{1+\tau^2} \cdot  Z ).
\]
But $\E \psi'_{\lambda/a} (Z) = B(\lambda/a,\Phi)$; so
\[
 \E \Psi'(\sqrt{1+\tau^2} \cdot  Z ; r)  = \left( \frac{r}{1+r} \right) \cdot B(\lambda/a,\Phi) .
\] 
We have the upper bound
$A(\lambda \cdot (1+r), H \star \Phi_\tau) \leq (\lambda \cdot (1+r))^2$ because $\|\psi_\kappa\|_\infty = \kappa$, and the
 lower bound
$B(\lambda \cdot (1+r), H \star \Phi_\tau) \geq 0$ because $\psi_\kappa' \geq 0$.
Moreover, both bounds are tight, as can be seen by choosing the point mass with
$H = \delta_\mu$ as $\mu \goto \infty$.
Combining all the above, we obtain the following. 
% and noticing that $\lambda/a = kappa (1+r)$

\begin{lemma} \label{lem:ABounds} With $r > 0$, $a \equiv \sqrt{1+\tau^2}/(1+r)$, and $F = (1-\eps) \Phi + \eps H$,
\begin{eqnarray*}
  \cA(\tau^2,r; m,\kappa,F) &=& (1-\eps) \left(  ar \right)^2 A(\kappa \cdot (1+r),\Phi).  + \eps \cdot  \left( \frac{r}{1+r} \right)^2 A(\lambda \cdot (1+r), H \star \Phi_\tau),\\
 \cA( \tau^2,r; m,\kappa, F)  &\leq&   (a r)^2 \cdot  \bar{A}(  \kappa \cdot  (1+r) , \eps ), \\
 \cA( \tau^2,r; m,\kappa, F) &\goto&   (a r)^2 \cdot  \bar{A}(  \kappa \cdot  (1+r) , \eps ), \qquad H = \delta_\mu, \quad \mu \goto \infty. 
 \end{eqnarray*}
\end{lemma}

\begin{lemma} \label{lem:BBounds} With $\cB = \E \Psi'( W +\tau Z ; r)$ and $ W \sim F = (1-\eps) \Phi + \eps H$,
\begin{eqnarray*}
  \cB(\tau^2; \dof,\kappa,F) &=& \left( \frac{r}{1+r} \right) \cdot
                                 \left ( (1-\eps) B(\kappa \cdot
                                 (1+r),\Phi) + \eps \cdot \E
                                 \psi'_{\lambda (1+r)} (U + \tau Z )
                                 \right )      \, , \\
  \cB( \tau^2,r; \dof,\kappa, F) &\geq&   \left( \frac{r}{1+r} \right)
                                        \cdot \bar{B} (\kappa \cdot
                                        (1+r),\eps)  \, ,\\
  \cB( \tau^2,r; \dof,\kappa, F) &\goto&   \left( \frac{r}{1+r} \right) \cdot \bar{B} (\kappa \cdot (1+r),\eps),  \qquad H = \delta_\mu, \quad \mu \goto \infty. 
\end{eqnarray*}
\end{lemma}

The proof of Lemma \ref{lem:extremalSE},
in the Appendix, combines the last two lemmas to 
obtain the equivalence of LFSE and double-bar evolution.

\section{Minimax Asymptotic Variance of Floating Threshold SE}

\subsection{Minimax Formal Variance}

\begin{definition}
Define the {\em formal  variance} 
\[
          \cV_\dof (\kappa,F) = \dof \cdot \tau_\infty^2(\dof,\kappa, F).
\]
where $\tau_\infty^2(\dof,\kappa, F)$ denotes the fixed point of the associated 
floating-threshold State Evolution.

Define the minimax formal variance to be
\[
          \cV^*_m(\eps) = \inf_\kappa \sup_{F \in \cF_\eps} \cV_\dof(\kappa,F).
\]
\end{definition}

The minimax problem identifies  
a distinguished  choice 
of the capping parameter, offering the best guarantee 
applicable across all $F \in \cF_\eps$.
Here is the solution:

\begin{lemma} \label{lem:calibrateKappa}
The mapping $\kappa \mapsto \barr(\kappa ; \dof,\eps)$ is continuous and strictly monotone decreasing.
For each $\bkappa > 0$,  the equation
\[
   \bkappa= \kappa \cdot (1 + \barr(\kappa))
\]
has an unique solution $\kappa = {\ukappa}(\bkappa)$.
\end{lemma}

\begin{thm} \label{thm:minmaxKappa}
Let $\kappa^*(\eps)$ denote Huber's minimax capping parameter
in the scalar estimation problem \cite{HuberMinimax}.
Let $\ukappa(\cdot)$ denote the re-calibrated function defined by
Lemma \ref{lem:calibrateKappa}.
Define the re-calibrated parameter
\[
\ukappa^*(\eps) ={\ukappa}(\kappa^*(\eps)).
\]

Suppose that $\dof \cdot I(F_\eps^*) > 1$; then {\em every}
 instance of floating-threshold state evolution having  parameters $(\dof,\tau_0,\ukappa^*(\eps) ,{F})$
 with proper $F \in \cF_\eps$ has a fixed point at $\tau_\infty^2 \equiv \tau_\infty^2(\dof,\ukappa^*(\eps) ,{F}) $
 obeying
 \[
      \tau_\infty^2 \leq   \bar{\tau}^2_\infty(m,\eps,\ukappa^*(\eps)) \equiv \frac{ v^*(\eps)/\dof}{1-\ v^*(\eps)/\dof}.
 \]
More specifically, we have the saddlepoint relation:
\[
     \inf_\kappa \sup_{F \in \cF_\eps} \cV_m(\kappa,F) = \cV_m (\ukappa^*(\eps),\bar{F}_\eps) =   \sup_{F \in \cF_\eps} \inf_\kappa \cV_m(\kappa,F),
\]
with saddle point at $(\ukappa^*(\eps),\bar{F}_\eps)$, and
where the minimax value $  \cV_m^*(\eps) = \cV_m (\ukappa^*(\eps),\bar{F}_\eps) $ obeys:
\[
    \cV_m^*(\eps) =  \frac{ v^*(\eps)}{1-\ v^*(\eps)/\dof}  \equiv \frac{1}{ I(F_\eps^*) -1/m}.
\]
\end{thm}

\begin{figure}[h]
\begin{center}
\includegraphics[height=3in]{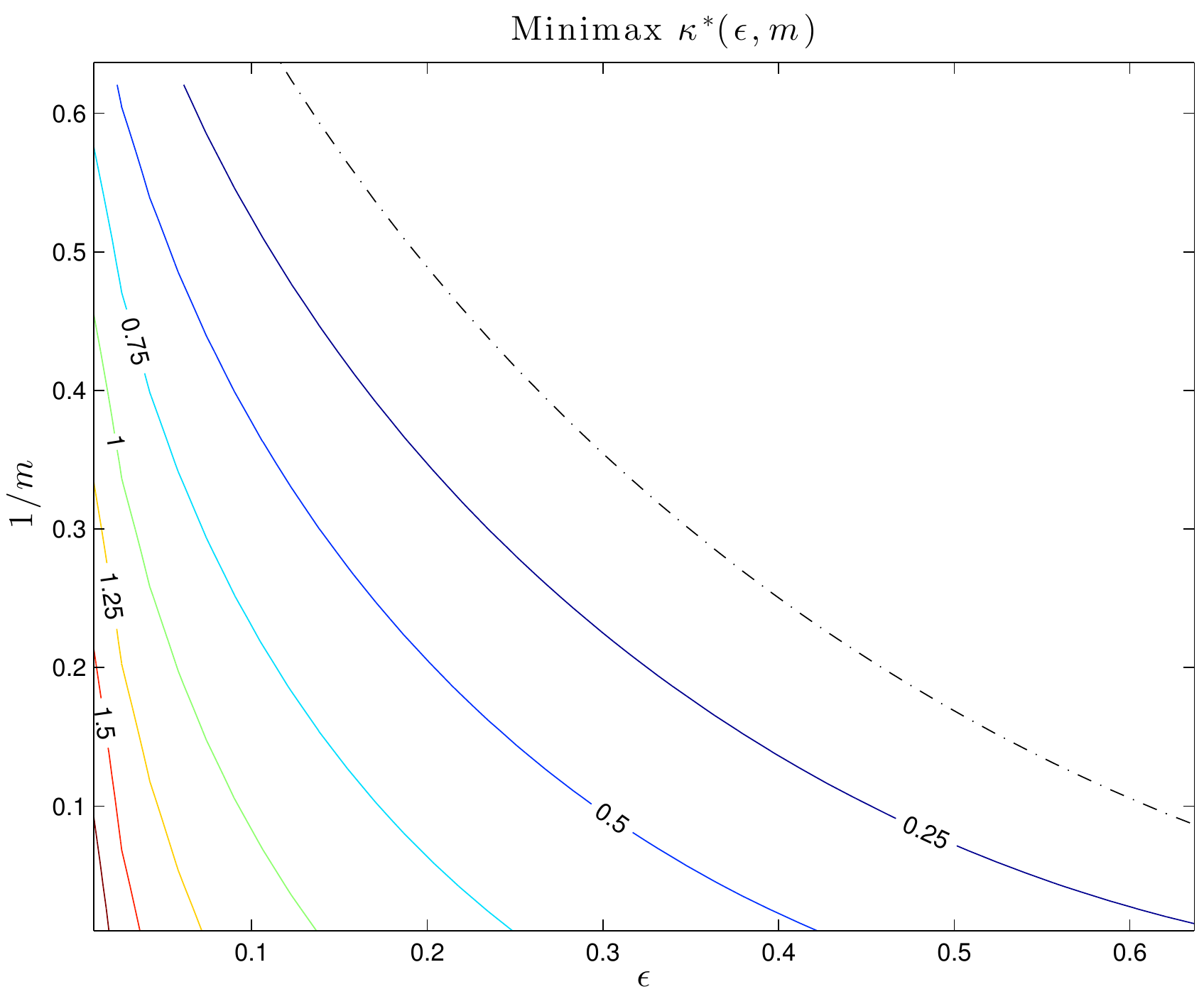}
\caption{Minimax $\ukappa^*(\eps;m)$.
Each pair $(\eps,m)$  is represented by the point
$x = \eps$ and $y = 1/m$. 
 Contours of the minimax capping parameter
$\ukappa^*(\eps;m)$ are depicted in the
region below the dashdot curve
at $1/m = I(F_\eps^*)$. }
\label{fig-Kappa}
\end{center}

\end{figure}

Figure \ref{fig-Fig2} presents a diagram
showing contours of $\cV_m^*(\eps)$.
The diagram employs  the unit square $\{ (\eps,1/m): 0 \leq \eps, 1/m \leq 1 \}$
where the $x$-axis shows the  contamination fraction $\eps$,
and the $y$ axis shows $1/m$ for plotting purposes.   Only the part of the diagram where $1/m < I(F_\eps^*)$
is populated with contours. The reader can see how the 
asymptotic variance `blows up' as $1/m$ approaches $I(F_\eps^*)$,

Figure \ref{fig-Kappa} 
shows contours of 
the minimax capping parameter $\ukappa^*(\eps;m)$.
The reader can see how the 
capping parameter  shrinks to zero as $1/m$ approaches $I(F_\eps^*)$,

\subsection{State Evolution in the Unbounded Phase}

Figure \ref{fig-Fig2} has a `bounded' phase, where the formal variance is
bounded across all contaminating distributions, 
and a complementary so-far undescribed phase. It seems that the formal
variance must be unbounded in this phase, since the phase consist of cases with smaller $m$ than the bounded ones,
and so therefore of `harder' cases. Validating this intuition, we have:

\begin{coro}  \label{lem:unboundedPhase}
Suppose that $\dof  \cdot I(F_\eps^*) \leq 1$; then for each $\tau < \infty$,  {and each $\kappa > 0$}, {\em some} 
 instance of proper state evolution with  parameters $(\dof,\tau_0,\kappa ,{F})$
 and proper $F \in \cF_\eps$ has a unique fixed point at $\tau_\infty^2 \equiv \tau_\infty^2(\dof,\kappa ,{F}) $
 obeying
 \[
      \tau_\infty^2  \geq    \tau^2.
 \]
\end{coro}

\begin{figure}[h]
\begin{center}
\includegraphics[height=3in]{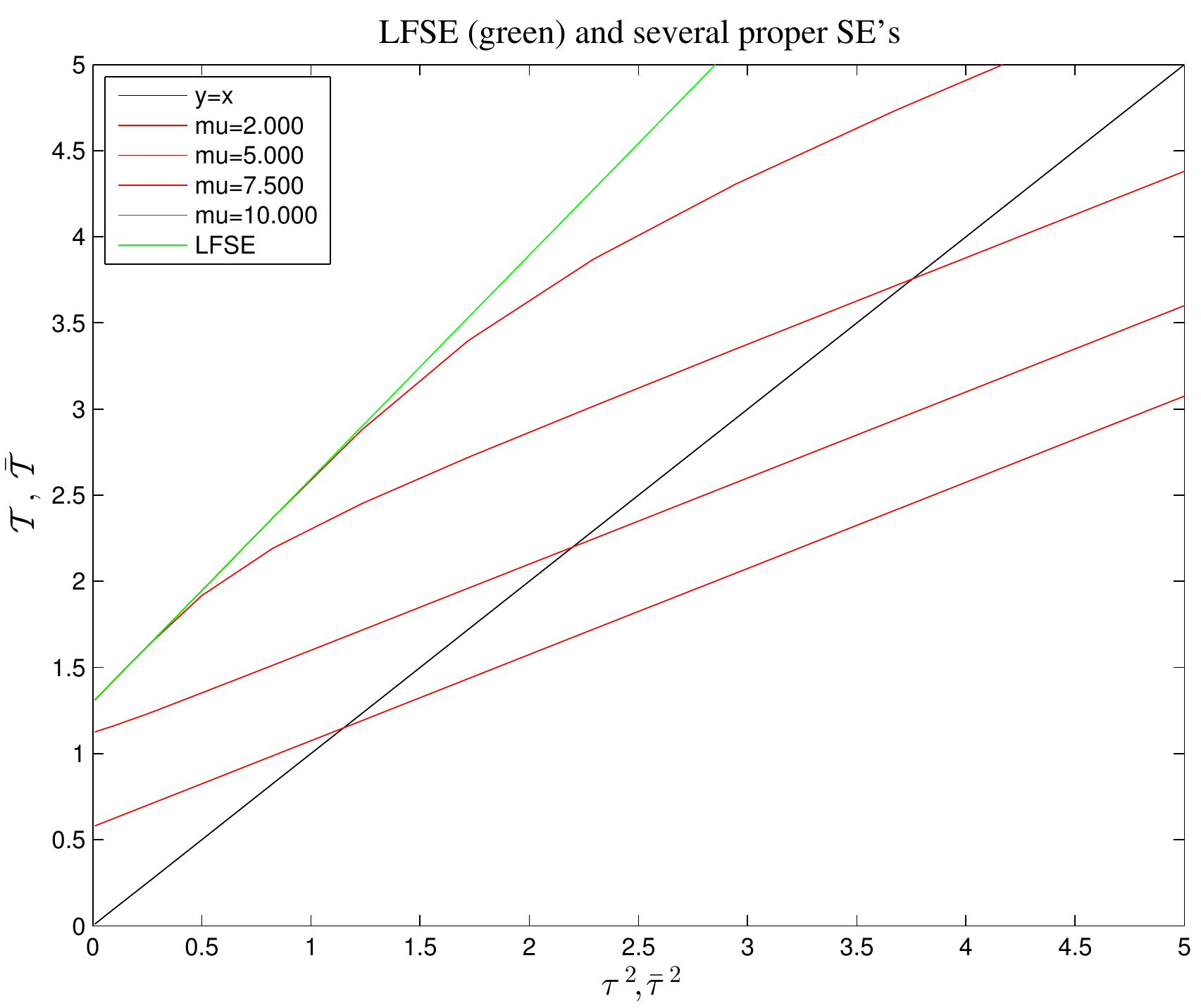}
\caption{SE in the unbounded phase. Examples of proper state evolutions with
$\mu = 2,5,7.5,10$, $\eps = 0.05$, $m =5$. The LFSE dynamical system has no fixed point.
The proper SE's have fixed points, but the location of the fixed point is unbounded above.}
\label{fig-Fig3}
\end{center}
\end{figure}

Goings-on in the unbounded phase are documented in Figure \ref{fig-Fig3}.
In the unbounded phase, {\em every} LFSE map $\bcT$
has no fixed point, whatever be the parameter $\kappa$. 
Proper state evolutions still have unique stable fixed points,
but there is no upper bound on their size. Hence the worst-case  fixed point $\btau_\infty^2$
is infinite.

This is an instance of what Donoho and Huber \cite{donoho1983notion} called  {\em  breakdown of asymptotic variance}. 
Breakdown occurs, in the $(\eps,1/m)$ phase diagram,  where-ever $m I(F_\eps^*) \leq 1$, and the breakdown point is $m I(F_\eps^*) = 1$,
the dashdot curve in our figures.

Note that as
$m \goto \infty$, we converge to the classical case, where the asymptotic variance of
$(M)$-estimates does not break down. In the high-dimensional case $n/p \goto m \in (1,\infty)$,
the asymptotic variance does break down.

\section{Minimax Variance of the Huber $(M)$-estimates}

We now develop our main result about $(M)$-estimates.
The analysis in the last section concerns
floating-$\lambda$ state evolution; while Theorem \ref{thm:SECorrectHuber}
shows that fixed-$\lambda$ state evolution describes 
the asymptotic variance of the Huber $(M)$-estimate.
We show how to bridge this difference.

\subsection{Minimax Formal Variance}

\begin{definition} \label{def:calibrate} {\bf Calibration Relation.}
Suppose the proper floating threshold
state evolution with parameters $(\dof,\tau_0,\kappa, F)$ has a unique fixed point $\tau_\infty^2$.
We formally associate this to a Huber $(M)$-estimate 
in the linear model under asymptotic regime $PL(m)$
with parameter $\lambda$  satisfying
\[
                  \lambda = \kappa \cdot \sqrt{1 + \tau_\infty^2(m,\kappa,F)}.
\]
\end{definition}

We denote this correspondence  by $\lambda = \lambda_\infty(m,\kappa,F)$ and 
the inverse correspondence with $\kappa = \kappa_\infty(m,\lambda,F)$.

\begin{definition}
The {\em  formal asymptotic variance} of the Huber $(M)$-estimator under the $PL(\dof)$ asymptotic framework is
\[
          \cV^\circ_\dof (\lambda,F) \equiv \dof \cdot \tau_\infty^2(\dof,\kappa, F),
\]
where $\tau_\infty^2(\dof,\kappa, F)$ denotes the fixed point of the floating threshold state evolution with parameter $\kappa$
and where  $\lambda = \lambda_\infty(m,\kappa,F)$.
\end{definition}

Theorem \ref{thm:SECorrectHuber} shows that this formula is rigorously correct --
the Huber estimator with the specified parameter $\lambda$ indeed has almost surely an asymptotic variance 
and it is equal to the
formal asymptotic variance.

\begin{lemma}  \label{lem:monotoneLambda}
Let $\bcV_m(\kappa,\eps) = \sup_{F \in \cF_\eps} \cV_m^\circ(\kappa,F)$
denote the worst case formal variance, across the full $\eps$-neighborhood,
of the  floating-threshold state evolution fixed point
under capping parameter  $\kappa$.
Set 
\[
\kappa^{+}(m,\eps) =  \sup \{ \kappa:   \bcV_{m}(\kappa,\eps) < \infty \};
\]
there is $m_0(\eps) \in (1,\infty)$ so that,
for $m > m_0(\eps)$, we have $\bcV_m(\kappa,\eps) <  \infty$ throughout
$(0,\kappa^{+}(m,\eps))$.
Define
\[
\blambda(\kappa; m,\eps) = \sup_{F \in \cF_\eps}  \lambda_\infty(m,\kappa,F).
\]
For $m > m_0(\eps)$, the  mapping 
\[
\kappa \mapsto   \blambda(\kappa; m,\eps)  % \kappa \cdot \bar{\sigma}(\kappa;m,\eps)  % \frac{\kappa}{\sqrt{1 - \bar{V}(\bkappa(\kappa),\eps)/m}}
\] 
is strictly increasing for $0 < \kappa < \kappa^+(m,\eps)$.
\end{lemma}

\begin{figure}[h]
\begin{center}
\includegraphics[height=3in]{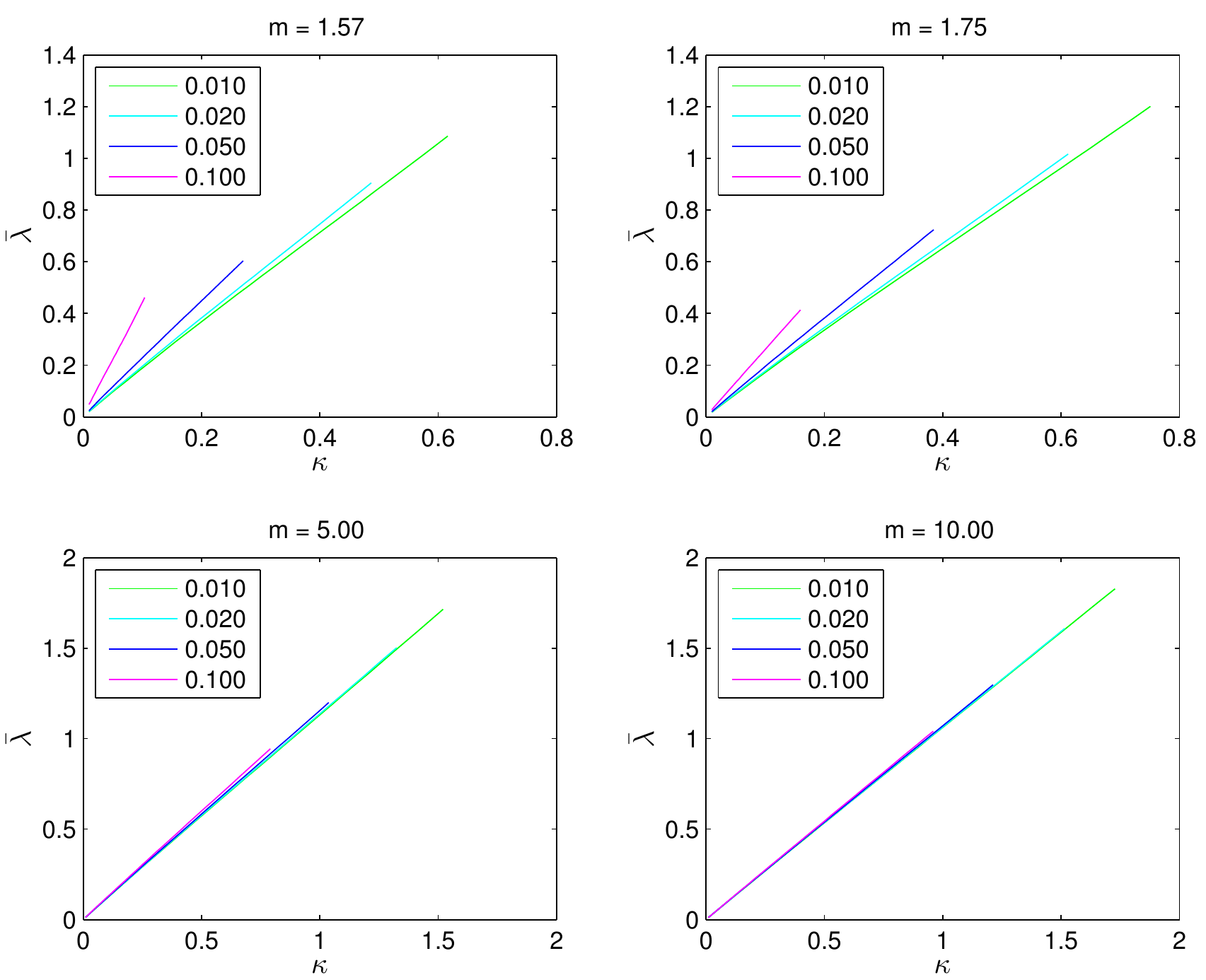}
\caption{Monotonicity of $\kappa \mapsto \bar{\lambda}(\kappa)$.
Each subplot depicts $\bar{\lambda}(\kappa;m,\eps)$
as a function of $\kappa$, for $\eps \in \{ 0.01,0.02,0.05,0.10 \}$,
at one particular $m$.  Evidently, as $m \goto \infty$,
$\bar{\lambda} \goto \kappa$.}
\label{fig:mono}
\end{center}
\end{figure}

Figure \ref{fig:mono}  displays $\blambda(\kappa)$ for a variety of choices of $\eps,m$;
the monotonicity is evident.
%Apparently $m_0$ is reasonably small in a wide range of cases.
Numerics show that we may take $m_0 \equiv \pi/2$;
however our proof only attempts to show that {\em some} $m_0$ sufficiently large will work.

The monotonicity condition on $\blambda$ ensures
that the least-favorable contamination for the Huber $(M)$-estimator is achieved by
the improper distribution $\bar{F}_{\eps}$.
\begin{thm}  \label{lem:minmaxThm} {\bf Evaluation of Minimax Asymptotic Variance of Huber (M)-estimator.}
If the mapping $\kappa \mapsto   \blambda(\kappa; m,\eps)$ is strictly increasing 
for $0 < \kappa < \kappa^+$ we have 
\[
       \inf_\lambda \sup_{F \in \cF_\eps} \cV^\circ_m(\lambda,F) = \inf_\kappa \sup_{F \in \cF_\eps} \cV_m^\circ(\kappa,F),
\]
where the minimax on the left concerns the formal variance of Huber $(M)$-estimates
parametrized by $\lambda$, and that on the right concerns the formal variance
of floating-threshold state evolutions parametrized by $\kappa$. The minimax tuning of the Huber (M)-estimator is achieved by the tuning parameter
\[
\lambda^*(\eps) = \blambda(m,\ukappa^*(\eps),\eps).
\]
\end{thm}
It follows of course that we have the formula
\[
  \inf_\lambda \sup_{F \in \cF_\eps} \cV^\circ_m(\lambda,F) = \frac{1}{I(F_\eps^*) -1/m},
\]
which agrees in the limit $m \goto \infty$ with Huber's classical formula
for the scalar location problem:
\[
  \inf_\lambda \sup_{F \in \cF_\eps} V(\psi_\lambda,F) = \frac{1}{I(F_\eps^*)}.
\]

Figure \ref{fig-Lambda} 
shows contours of 
the minimax thresholding parameter $\lambda^*(\eps;m)$.
The reader can see how this 
parameter  shrinks to zero as $1/m$ approaches $I(F_\eps^*)$.
While the story is much the same as for the $\kappa$ parameter
in Figure \ref{fig-Kappa}, the $\lambda$-parameter is the one relevant to
practice, because the $\kappa$ parameter is a theoretical construct
while the corresponding $\lambda$ parameter can actually be used to 
specify the desired Huber estimator in statistical software packages.

Since the formal variance
$\cV^\circ_\dof (\lambda,F)$ has the saddlepoint property,
Theorem 2.2 shows  that  the rigorous asymptotic variance
$\AVar = \AVar(\hat{\theta}^\lambda_n,F)$ (say) has it as well.
 
 \begin{definition}
Let $\cF_\eps^2$ denote the subset of distributions in $\cF_\eps$ with finite variance:
$\mu_2(F) = \int w^2 \de F(w) < \infty$.
\end{definition}

\begin{coro} \label{coro-minimax}
Fix $\eps > 0$, and $m > m_0(\eps)$. We are in the 
asymptotic regime $PL(m)$.
\bitem 
\item Suppose that $\cV_m^*(\eps)$ is finite.
Consider the formally minimax parameter  $\lambda = \lambda^*(\eps)$;
let $\hat{\theta}_n^*$ denote a corresponding solution 
of the Huber (M)-equation with that $\lambda$. 
For every  error distribution $F \in \cF^2_\eps$, we have
\[
  \AVar(\hat{\theta}_n^*,F)   \leq   \cV^*_m(\eps)\, .
\] 
For every  $\lambda \neq \lambda^*(\eps)$ there is a proper $\eps$-contaminated
normal error distribution $F \in \cF_\eps^2$ so that the Huber estimator $\hat{\theta}^{\lambda}$
obeys
\[
   \AVar (\hat{\theta}^{\lambda}, F)  > \cV^*_m(\eps)\, .
\] 
Consequently,
\[
    \inf_\lambda \sup_{F \in \cF_\eps^0} \AVar(\hat{\theta}^{\lambda},F) =
      \AVar( \hat{\theta}_n^*, \bar{F}_\eps) =  \cV^*_m(\eps) \, .
\]
\item Suppose that $\cV_m^*(\eps)$ is infinite.
For every  $\lambda > 0$ and each $V > 0$, 
there is a proper $\eps$-contaminated
normal error distribution $F \in \cF_\eps^2$ with 
\[
   \AVar(\hat{\theta}^{\lambda}, F)  >  V \, .
\] 
Consequently,
\[
    \inf_\lambda \sup_{F \in \cF_\eps^2} \AVar(\hat{\theta}^{\lambda},F) =
      +\infty =  \cV^*_m(\eps) .
\]

\eitem
\end{coro}

\section{Discussion}

Under  the high-dimensional $PL(\dof)$ asymptotic  \-- as shown in \cite{bean2013optimal} \-- 
the maximum likelihood estimator  is no longer an efficient
 estimator.  It follows that the Huber estimator is no longer
 asymptotically minimax among all (M)-estimators. Hence the asymptotic minimax
 in (\ref{eq:modern}) should better be called the asymptotic  minimax among Huber estmates.
The degree of sub optimality can be controlled explicitly.
By \cite[Corollary 3.7]{donoho2013high}, the asymptotic variance
under $PL(\dof)$ obeys the following inequality,
which is strictly stronger than the Cram\'{e}r-Rao bound when $ 1 < \dof < \infty$:
\[
V_\dof(\psi_\lambda, F)  \geq  \frac{1}{1-1/\dof}  \cdot  \frac{1}{I(F)},
\]
and so the minimax asymptotic variance obeys:
\beq \label{eq:suboptimal}
    \min_\lambda \max_{F \in \cF_\eps}  V_\dof(\psi_\lambda, F)  \geq  \frac{1}{1-1/\dof}   \cdot \frac{1}{I(F_\eps^*)} .
\eeq
It follows that provided $m I(F_\eps^*) > 1$, 
\[
\mbox{ minimax Huber asymptotic variance} \leq  K \cdot \mbox { minimax asymptotic variance} ,
\]
where
\[
   K = K(\dof,\eps) = \frac{1-1/\dof}{1-I(F_\eps^*)/\dof}.
\]
One sees directly that the sub-optimality of the Huber estimator is well controlled
provided that $I(F_\eps^*)$ is close to one; i.e., in the regime where $\eps$ is small enough (though this is
$m$-dependent).  Of course in the regime $I(F_\eps^*) \dof \leq 1$,
some other estimators could be dramatically more robust.

\section*{Acknowledgements}

This work was partially supported by:  NSF grants CCF-1319979 (A.M.);
DMS-1418362 and DMS-1407813 (D.D.);
and the grant
AFOSR FA9550-13-1-0036.

\section*{Appendix: Proofs}

\subsection*{Proof of Lemma \ref{lem:dominance}}
The desired relations are true for iteration $t=0$ by assumption (note that
no assertion about the sequence $(r_t)_{t \geq 1}$
 is made at stage $t=0$, only about $\tau_0$).

Suppose that we have proved the desired relations up to iteration $t-1$ and we now must show
that they hold for iteration $t$.

We observe  that $\bar{F}_{\eps,\sigma}$ is stochastically more spread than
any proper distribution in $\cF_{\eps,\sigma}$ \-- that is, every distribution with
all its mass on the reals $\reals$ rather than the extended reals $\reals \cup \{\pm \infty \}$.
Hence for every function $\xi(x)$ monotone increasing  in $|x|$,
\[
   \sup\{ \E  \xi(X) : X \sim F \in \cF_{\eps,\sigma} \}  \equiv \E  \xi ( \bar{X}_{\eps,\sigma}).
\]
In this sense $ \bar{X}_{\eps,\sigma}$ is extremal among contaminated normals.
Moreover, we note that for $\tsigma > \sigma$, $\bar{F}_{\eps,\tsigma}$
is more spread than  $\bar{F}_{\eps,{\sigma}}$.
Hence, again for $\xi$ monotone increasing in $|x|$,
\[
\E  \xi ( \bar{X}_{\eps,\sigma}) \leq \E  \xi ( \bar{X}_{\eps,\tsigma}).
\]

Applying Lemma \ref{lem:monotonicity}, Claim 2,
\beq \label{eq:PsiP}
 \inf \{ \E_{F} \Psi'_{\lambda} (X,r) : X \sim F \in \cF_{\eps,\sigma} \}  \equiv \E \Psi'_{\lambda} (\bar{X}_{\eps,\sigma} , r), 
   \quad \forall \, b, \lambda >0.
\eeq
So in particular, for $X \sim F \in \cF_{\eps,\sigma_{t-1}}$ 
we have $\E \Psi'(\bar{X}_{\eps,\sigma_{t-1}} , r_t) \leq \E \Psi'(X , r_t)$.
Hence we must have
\begin{eqnarray*}
   \frac{1}{\dof} &=&  \E \Psi'_{\lambda_t} (X, r_t) \\
   & \geq &  \E \Psi'_{\lambda_t} (\bar{X}_{\eps,\sigma_{t-1}}, r_t) \\
      & = &  \E \Psi'_{\tlambda_t} (\bar{X}_{\eps,\tsigma_{t-1}},
            r_t)\, .
\end{eqnarray*}
The first step is just the definition of $r_t$,
the second step  used (\ref{eq:PsiP}), and the third step used (\ref{eq:invar}).
Now since $r \mapsto \frac{r}{1+r}$ is increasing in $r$,  while   
$r \mapsto \P\{ |\bar{X}_{\eps,\tsigma_{t-1}}| >  \kappa \cdot \tsigma_{t-1} \cdot ( 1+r)  \}$ is monotone increasing in $r$. 
Hence the product  \-- $r \mapsto \E \Psi'_{\tlambda_t} (\bar{X}_{\eps,\tsigma_{t-1}}, r) $ \-- 
is monotone increasing in $r$; so in order to satisfy the definition of $\tarr_t$ \-- 
\[
 \frac{1}{\dof} = \E \Psi'_{\tlambda_t} (\bar{X}_{\eps,\tsigma_{t-1}}, \tarr_t) 
\]
 \--  we  must have
\[
\tarr_t \geq r_t.
\]

Now turn to the dominance relation concerning $\tau_t^2$.
\begin{eqnarray*}
 \sup\{ \E \Psi_{\lambda_t}^2 (X,r_t) : X \sim F \in \cF_{\eps,\sigma_{t-1}} \}  & \equiv & \E \Psi_{\lambda_t}^2 (\bar{X}_{\eps,{\sigma}_{t-1}},r_t) \\
 &=& \left( \frac{\sigma_{t-1}}{\tsigma_{t-1}} \right)^2 \cdot  \E \Psi_{\tlambda_t}^2(X_{\eps,\tsigma_{t-1}},r_t) \\
   &\leq&  \E \Psi_{\tlambda_t}^2 ( \bar{X}_{\eps,\tsigma_{t-1}} , \tarr_t).
\end{eqnarray*}
where in the first inequality we substituted (\ref{eq:quad-invar})
and in the second inequality we substituted 
$r_t \mapsto \tarr_t$ by Lemma \ref{lem:monotonicity}, Claim 3. We conclude that
\[
   \tau_t \leq \ttau_t,
\]
which completes iteration $t$ of the claimed result and sets up the
assumptions for the next iteration. \qed

\subsection*{Proof of Lemma \ref{lem:monotonicity}}

To prove the $i$-th claim, for $i=1,\dots, 4$, combine
formula  $\Psi(x;r) = r\, \psi_\lambda(x/(1+r))$
with the corresponding numbered observation:
\begin{enumerate}
\item From $  |\psi_\lambda(x)| = \min(|x|,\lambda)$, $|x| \mapsto |\psi_\lambda(x)|$ is nondecreasing.
\item From $ \psi'_\lambda(x) = 1_{\{ |x| \leq \lambda \}}$, $|x| \mapsto  \psi'_\lambda(x) $ is nonincreasing.
\item $r \mapsto r/(1+r)$ is monotone increasing.
\item $\lambda \mapsto \min(|x|,\lambda)$ is nondecreasing.
\end{enumerate}
\qed

\subsection*{Proof of Lemma \ref{lem:invariances}}

These are simple scaling invariances, combined with $\Psi(x;r) = \frac{r}{1+r} \psi_\lambda(x)$.

The first, (\ref{eq:invar}), says simply that for any $0 < \sigma < \tsigma$,
\begin{eqnarray*}
    \E \psi'_{\kappa \sigma } (\bar{X}_{\eps,\sigma}) &=&  \P\{ |\bar{X}_{\eps,\sigma} | < \kappa \cdot \sigma \}  \\
    &=& \Phi_\sigma(-\kappa \sigma, \kappa \sigma) \\
    &=& \Phi(-\kappa,\kappa) \\
    &=&    \Phi_{\tsigma}(-\kappa \tsigma, \kappa \tsigma)  \\
    &=&   \P\{ |\bar{X}_{\eps,\tsigma} | < \kappa \cdot \tsigma \} \\
&=& \E \psi'_{\kappa \tsigma } (\bar{X}_{\eps,\tsigma}) .
\end{eqnarray*}

The second, (\ref{eq:quad-invar}),  combines two invariances. If $Z \sim N(0,1)$, then
\[
    \E \psi^2_{\kappa \sigma } (\sigma Z) =  \sigma^2 \cdot \E \psi^2_{\kappa  } (Z)  ;
\]
while, if $U \sim H_\infty$ is the degenerate improper random variable supported at $\infty$,
\[
   \E \psi^2_{\kappa \sigma } (U) = \kappa^2 \cdot \sigma^2  .
\]
Hence 
\[
    \E \psi^2_{\kappa \sigma } (X_{\eps,\sigma}) = \sigma^2 \cdot ( (1-\eps) \E \psi^2_{\kappa  } (Z)  + \eps \kappa^2) ,
\]
is thus proportional to $\sigma^2$.
Applying this both to $\sigma = \sigma_{t-1}$ and $\sigma = \tsigma_{t-1}$ gives (\ref{eq:quad-invar}).
\qed

\subsection*{Proof of Lemma \ref{lem:extremalSE}.}

Note that $r \mapsto \frac{r}{1+r}$ and $r \mapsto B(\kappa \cdot (1+r), F)$ 
are each strictly monotone increasing in $r > 0$. Moreover, 
\[
  \inf_{F \in \cF_{\eps} } B(\kappa \cdot (1+r), F) = \bar{B} (\kappa \cdot (1+r),\eps)  .
\]
Set $R(b) = 1/(m b -1)$; this is monotone decreasing in $b >1/m$.
Then  $\cR =  R(\cB)$ and  this relationship  is monotone decreasing  in $\cB > 1/m$.
Hence
\[
    \sup_{F \in \cF_{\eps} } R( \cB(\tau^2,r;m\eps,\kappa)) = R(\bar{B} (\kappa \cdot (1+r),\eps)
 \]
Also  note that $\bbarr=  R(\bar{B} (\kappa \cdot (1+\bbarr),\eps)$.
It follows that for each $\eta>0$, 
for some $r > \bbarr - \eta$ we can find $F \in \cF_\eps$ satisfying
\[
 r =  R( \cB(\tau^2,r;m,\kappa,F)).
\]
Hence,
\[
\sup_{F \in \cF_{\eps} } \cR(\tau^2;m,\kappa,F)\geq \bbarr.
\]
Since $r \mapsto \bar{B}(\kappa (1+r),\eps)$ is strictly monotone increasing, we conclude
that for  every $r > \bbarr$ we must have
\[
    r > \bbarr=  R(\bar{B} (\kappa \cdot (1+\bbarr),\eps)) > R(\bar{B} (\kappa \cdot (1+r),\eps)) 
\]
and hence, for such $r$
\[
    r >  \sup_{F \in \cF_{\eps} } R( \cB(\tau^2,r;m,\kappa,F)), 
\]
implying that 
\[
r > \sup_{F \in \cF_{\eps} } \cR(\tau^2;m,\kappa,F),
\]
and so also 
\[
\bbarr \geq \sup_{F \in \cF_{\eps} } \cR(\tau^2;m,\kappa,F),
\]
which proves (\ref{eq:max-r}).

We turn to Eq. (\ref{eq:max-A}).
Set $\bar{\bar{a}} \equiv \sqrt{1+\tau^2}/(1+\bbarr)$. We have
\begin{eqnarray*}
 \sup_{F \in \cF_{\eps} } \cA( \tau^2, \cR(\tau;m,\kappa,F ) ; m,\kappa, F)   &=&  
  \sup_{F \in \cF_{\eps} } \cA( \tau^2,\bbarr ; m,\kappa, F)  \\
 &=&  ( \bar{\bar{a}} \bbarr)^2 \cdot  \bar{A}(  \kappa \cdot  (1+\bbarr) , \eps ) \\
  &=&  (1 + \tau^2) \left( \frac{\bbarr}{1+\bbarr} \right)^2 \bar{A}(  \bkappa , \eps ) \\
  &=& (1 + \tau^2)  \frac{ \bar{A}(  \bkappa , \eps )}{ m^2 \bar{B}(  \bkappa, \eps )} \\
  &=& (1 + \tau^2) \cdot \bar{V}( \bkappa,\eps) /m^2 \\
   &\equiv& \bbcA(\tau^2;m,\eps,\kappa).Œ
\end{eqnarray*}
In the first step we used monotonicity of 
$r \mapsto \cA(\tau^2, r; m,\kappa,F)$, 
and  in the second step,
we used Lemma \ref{lem:ABounds}. In each step inequality is clear, while
equality is demonstrated by choosing a  sequence of contamination cdfs
 $H = \delta_\mu, \quad \mu \goto \infty$. \qed

\subsection*{Proof of Lemma \ref{lem:calibrateKappa}}

For fixed $(\eps,m)$, consider the  relationship between $(\bbarr,\bkappa)$ implied by
\[
  \frac{\bbarr}{1+\bbarr} \bar{B}(\eps,\bkappa) = \frac{1}{m} .
\]
Note that $ \bar{B}(\eps,\bkappa) = (1-\eps) (2 \Phi(\bkappa) -1)$ where $\Phi(x)$ is the standard normal CDF,
which is a bijection between $(-\infty,\infty)$ and $(0,1)$.
One can check that, for fixed $(\eps,m)$, $(\bbarr, \bkappa)$ are in one-one correspondence
by the functions
\[
    \bbarr(\bkappa) = \frac{(1-\eps) (2 \Phi(\bkappa) -1)}{1/m- (1-\eps) (2 \Phi(\bkappa) -1)}
\]
and
\[
   \bkappa(\bbarr) = \Phi^{-1} \left(  ( 1 + \frac{1 + 1/\bbarr}{m (1-\eps)} ) /2 \right), 
\]
acting as bijections $\bbarr \leftrightarrow \bkappa $ between domains $(0,\infty)$
and $(0,\bkappa^*)$, where $\bkappa^*(\eps,m) = \Phi^{-1}( ( 1 + \frac{1}{m (1-\eps)} ) /2)$.
Defining
\[
   \kappa(\bkappa) = \bkappa / (1 + \bbarr(\bkappa)),
\]
the pair $(\kappa(\bkappa),\bbarr(\bkappa))$  will obey the relation  
\[
 \frac{\bbarr}{1+\bbarr} \bar{B}(\eps,\kappa(1 + \bbarr)) = \frac{1}{m} .
\]
We obtain the explicit expression
\[
   \ukappa(\bkappa) = \frac{\bkappa}{ 1+ \bbarr(\bkappa) } ,
\]
showing directly that $\kappa$ is uniquely defined in terms of
$\bkappa$, for given $(\eps,m)$.
\qed

%Recall $\barr = R(B(\kappa(1+\barr),\eps))$
%and notice that
%$\kappa \mapsto B(\kappa(1+r),\eps)$ is strictly increasing.
%For $x > 1/m$, $R(x) = 1/(m \cdot x-1)$ obeys $R'(x) < 0$.
%It follows that 
%\[
%     \frac{d}{d\kappa}  \bar{r} < 0.
%\]
%
%{\bf COMPLETE THIS LEMMA}

\subsection*{Proof of Lemma \ref{thm:minmaxKappa}}

By Lemma \ref{lem:extremalSE},
the variance map $\bcT$ is the pointwise supremum of all
variance maps of proper floating-threshold state evolutions with $F \in \cF_\eps$.
Hence, no proper FTSE can have a 
larger fixed point; i.e.
\[
   \tau_\infty^2(m,\kappa,F) \leq \btau_\infty^2(m,\kappa,\eps), \qquad \forall F \in \cF_\eps.
\]
From $\cV_m(\kappa,F) = m \cdot   \tau_\infty^2(m,\kappa,F)$, we have
\[
  \sup_{F \in \cF_\eps} \cV_m(\kappa,F)    =  m \cdot \btau_\infty^2(m,\kappa,\eps) ,
\]
and so, if $\btau_\infty^2(m,\kappa,\eps)  <\infty$  \-- implying  $\bar{V}(\bkappa(\kappa),\eps) < m$ \-- 
\[
\sup_{F \in \cF_\eps} \cV_m(\kappa,F)  = \frac{\bar{V}(\bkappa(\kappa),\eps)}{1 - \bar{V}(\bkappa(\kappa),\eps)/m} .
\]
Setting $\cK_0 = \{ \kappa :  {\bar{V}(\bkappa(\kappa),\eps)} < m \}$,
\[
\inf_\kappa \sup_{F \in \cF_\eps} \cV_m(\kappa,F)  =  \inf_{\kappa \in \cK_0}  \frac{\bar{V}(\bkappa(\kappa),\eps)}{1 - \bar{V}(\bkappa(\kappa),\eps)/m} .
\]
Now by construction, 
\begin{equation} \label{eq:optim}
\bar{V}(\bkappa(\ukappa^*(\eps)),\eps) = \bar{V}(\kappa^*(\eps),\eps) = v^*(\eps);
\end{equation}
and moreover for $\kappa \neq \ukappa^*(\eps)$, 
$\bkappa(\kappa) \neq \kappa^*(\eps)$; so 
\[
\bar{V}(\bkappa(\kappa),\eps) = \bar{V}(\bkappa(\kappa),\eps) > \bar{V}(\kappa^*(\eps),\eps) = v^*(\eps).
\]
Now $v \mapsto v/(1-v/m)$ is monotone increasing on $\{v: v\leq m\}$.
Consequently, if $v^*(\eps) < m$
\[
\inf_\kappa \sup_{F \in \cF_\eps} \cV_m(\kappa,F)  =   \frac{v^*(\eps)}{ 1 - v^*(\eps)/m}.
\]
By hypothesis $v^*(\eps)/m  \equiv 1/(m I(F_\eps^*)) < 1$,  and so this formula indeed holds.

Now note that automatically
\[
\inf_\kappa \sup_{F \in \cF_\eps} \cV_m(\kappa,F) \geq  \sup_{F \in \cF_\eps} \inf_\kappa  \cV_m(\kappa,F)  ;
\]
hence the argument will be completed by showing that
\[
  \sup_{F \in \cF_\eps} \inf_\kappa  \cV_m(\kappa,F)   =  \frac{v^*(\eps)}{ 1 - v^*(\eps)/m}.
\]

But we have already shown by (\ref{eq:optim}) that
\[
  \inf_\kappa  \cV_m(\kappa,\bar{F}_\eps) =  \frac{v^*(\eps)}{ 1 - v^*(\eps)/m} .
\]
For all but purists, this completes the proof of the saddlepoint relation
\[
  \sup_{F \in \cF_\eps} \inf_\kappa  \cV_m(\kappa,F)   =  \frac{v^*(\eps)}{ 1 - v^*(\eps)/m} =   \inf_\kappa  \cV_m(\kappa,\bar{F}_\eps) .
\]

Purists who want everything stated using proper RV's
will want the following spelled out.
Let $G_{\eps,\mu} = (1-\eps) \Phi + \eps H_\mu$.
 For $\eta > 0$, there is $\mu \in \reals$
with
\[
  \inf_\kappa  \cV_m(\kappa,G_{\eps,\mu}) > \inf_\kappa  \cV_m(\kappa,\bar{F}_\eps)  - \eta .
\]
Now note also that, for $\mu >\kappa$, the Huber $\psi_\kappa$ obeys
\[
     V(\kappa, G_{\eps,\mu}) =   V(\kappa, G_{\eps,\infty}  ) \equiv V(\kappa, \bar{F}_\eps ) = \bar{V}(\kappa,\eps) > v^*(\eps) .
\]
with similar statements also being true for $A$ and $B$.
This observation can be elaborated into a full proof, exploiting
\[
  V(\kappa, (1-\eps) \Phi + \eps N(\mu,\gamma) ) \goto V(\kappa, \bar{F}_\eps )
\]
as $\mu \goto \infty$. We omit the details. \qed

\subsection*{Proof of Lemma \ref{lem:unboundedPhase}}

If $m \cdot I(F_\eps^*) \leq 1$, then $\bar{V}(\kappa,\eps)/m \geq v^*(\eps)/m = 1/(m I(F_\eps^*)) \geq 1$
and so $ \bar{V}(\bkappa(\kappa),\eps)/m \geq 1$ for each $\kappa > 0$.

The  variance map of the  LFSE with parameters $(m,\kappa,\eps)$
is affine: 
\[
\bcT(\tau^2)  = \frac{\bar{V}(\bkappa(\kappa),\eps)}{m} (1 + \tau^2);
\]
so both the slope and intercept equal $\bar{V}(\psi_{\bkappa},\eps)/m \geq 1$.
Hence there is no fixed point, and in fact there is a strict vertical gap between 
the identity line and the graph of $\bcT$  \-- a gap of size $\bar{V}(\psi_{\bkappa},\eps)/m \geq 1$.

Now $\bcT$ is the pointwise supremum of all the variance maps of proper state evolutions.
Hence for any $\tau$ we choose, there is a variance map $\cT$ of some proper state
evolution lying above the diagonal line at $\tau^2$:
\[
  \cT(\tau^2) > \tau^2,
\]
which implies that the corresponding
highest fixed point $\cT(\tau_\infty^2) = \tau_\infty^2$ obeys $\tau_\infty^2 > \tau^2$.
For all but purists, this completes the proof.

Purists will want to know that among the highest such fixed points  are
 in fact unique fixed points, which then represent variances that are in fact achieved.
 We will show this for contaminated distributions of the form
 $G_{\eps,\mu}$, for large $\mu$.
 
 %Combining these observations with (\ref{eq:rescale-A})-(\ref{eq:rescale-B}),
 %the definitions of $\cA$ and $\cR$, and the scaling (\ref{eq:scaleInvar}), we conclude that 

 \bitem
 \item For such $G_{\eps,\mu}$,  $\mu$ sufficiently large,
 we will show that the variance map $\cT$ is
 star shaped; namely, defining $T = T(\tau^2 ; G_{\eps,\mu})$ by
 \[
       \cT(\tau^2) = (1 + \tau^2) \cdot T(\tau^2),
 \]
 then we will show that for 
 $\mu$ large, $\tau^2 \mapsto T(\tau^2)$ 
 is a monotone nonincreasing function of $\tau^2$.

 \item Any such star-shaped map has a unique fixed point; if $\tau_1^2 < \tau_2^2$ are two distinct
 purported fixed points then because the line 
 $\tau^2 \mapsto (1+\tau^2) T(\tau_1^2)$ has a unique fixed point at $\tau^2 = \tau_1^2$,
 then
\begin{equation} \label{eq:fpfree}
    \tau^2  >  (1+\tau^2) T(\tau_1^2), \qquad \tau^2 > \tau_1^2 .
\end{equation}
 Hence
 \begin{eqnarray*}
 (1 + \tau_2^2) T(\tau_2^2) & \leq &  (1 + \tau_2^2) T(\tau_1^2) \\
    & = & (1 + \tau^2) T(\tau_1^2) |_{\tau^2 = \tau_2^2} \\
    & < &  \tau_2^2 .
 \end{eqnarray*}
 In the last step we use (\ref{eq:fpfree}), evaluated at $\tau^2 = \tau_2^2$.
The last display contradicts the supposed fixed-point nature
of $\tau_2^2$ and  proves that the second fixed point
$\tau_2^2$ cannot exist.

 \eitem

To explain the star-shapedness,
we need to develop some rescaling relationships.
Let $S^\sigma F$ denote the rescaling operator on CDF's,
producing $(S^\sigma F)(x) = F(x/\sigma)$. 
For a given $F \in \cF_{\eps}$ and a given
$\tau^2$ and associated $\sigma^2 = 1+\tau^2$,
let $\tilde{F}^\sigma \equiv S^\sigma (F \star \Phi_\tau)$. 
We then have
\[
  \tilde{F}^\sigma = (1-\eps) \Phi + \eps \tilde{H}^\sigma ,
\]
where the contamination CDF $\tilde{H}^\sigma = S^\sigma(H \star \Phi_\tau)$.
Because of the scale invariance $\lambda = \kappa \sigma$,
\beq \label{eq:scaleInvar}
     \int \psi_{\kappa \sigma}^2 ( x ) \de (F \star \Phi_\tau)(x) = \sigma^2 \int \psi_{\kappa} ( u ) \de \tilde{F}^\sigma(u).
\eeq
Similarly,
\[
     \int \psi_{\kappa \sigma}' ( x ) \de (F \star \Phi_\tau)(x) =  \int \psi_{\kappa}'( u ) \de \tilde{F}^\sigma(u).
\]

It follows that
 \[
 T(\tau^2) \equiv  \frac{T(\tau^2)}{1+\tau^2}  \equiv  \frac{T(\tau^2)}{\sigma^2}  = (\frac{r}{1+r})^2 \cdot A(\psi_{\tkappa}, \tilde{F}^\sigma),
 \]
 where $\tkappa = \kappa(1+r)$ solves
 \[
   \frac{r}{1+r} \cdot B(\psi_{\tkappa}, \tilde{F}^\sigma) = \frac{1}{m} .
 \]

The reader should check that
the following claims, if established,
would combine to prove the desired monotonicity of $T$.
\bitem
\item $A(\psi_{\kappa_0},\tilde{F}^\sigma)$ is monotone decreasing in $\sigma$, for fixed $\kappa_0$.
\item $B(\psi_{\kappa_0},\tilde{F}^\sigma)$ is monotone decreasing in $\sigma$, for fixed $\kappa_0$.
\item $\kappa \mapsto A(\psi_\kappa,F)$ is increasing in $\kappa$.
\item $\kappa \mapsto B(\psi_\kappa,F)$ is increasing in $\kappa$.
\item $\sigma \mapsto r$ is monotone decreasing in $\sigma$.
\item $\sigma \mapsto \tkappa$ is decreasing in $\sigma$.
\eitem
Some of these are obvious  \-- for example, 
 monotonicity of $\kappa \mapsto A(\psi_\kappa,F)$  and $\kappa \mapsto B(\psi_\kappa,F)$.
 Others follow from earlier items  \-- monotonicity of $\sigma \mapsto \tkappa$ follows from
 that of $\sigma \mapsto r$, while monotonicity of $\sigma \mapsto r$ follows from
 the two earlier claims about $B$. Finally, the first two claims will be shown for $F = G_{\eps,\mu}$ 
 for all sufficiently large $\mu$. 
 
In the coming two paragraphs, let $\kappa$ be fixed independent of $\sigma$.
Now of course
\begin{eqnarray}
 A(\psi_{\kappa}, \tilde{F}^\sigma) &=&     \int \psi_{\kappa}^2 ( u ) \de \tilde{F}^\sigma(u) \\
 &=& (1-\eps) \int  \psi_{\kappa}^2(u) \de \Phi 
      + \eps \int  \psi_{\kappa}^2(u) \de \tilde{H}^\sigma(u);  \label{eq:rescale-A} \\
&=& I + II. \nonumber
\end{eqnarray}
the term $I$ being independent of $\sigma$, we focus on the second one, $II$. Similarly,
\begin{eqnarray}
   B(\psi_{\kappa}, \tilde{F}^\sigma)  
   &=&  \int \psi_{\kappa}' ( u ) \de \tilde{F}^\sigma(u) \\
   &=& (1-\eps) \int  \psi_{\kappa}'(u) \de \Phi 
      + \eps \int  \psi_{\kappa}'(u) \de \tilde{H}^\sigma(u);  \label{eq:rescale-B}\\
&=& III + IV. \nonumber
\end{eqnarray}
We again focus on the $\sigma$-varying term; this time $IV$.
Letting $H^{\sigma} = S^\sigma H$ we have
$\tilde{H}^\sigma= \Phi_{\tau/\sigma} \star H^\sigma$.
By associativity of convolution,
\[
     \int \psi_{\kappa}^2 ( u ) \de \tilde{H}^\sigma(u) =  \int
     (\psi_{\kappa}^2 \star \Phi_{\tau/\sigma}) ( u ) \de
     H^\sigma(u)\, .
\]
Similarly,
\[
     \int \psi_{\kappa}' ( x )   \de \tilde{H}^\sigma(u)=  \int (\psi_{\kappa}'\star \Phi_{\tau/\sigma}) ( u ) \de {H}^\sigma(u).
\]

Now note that, for all sufficiently large $u$, $u \mapsto (\psi_{\kappa}^2 \star \Phi_{\tau/\sigma}) ( u ) $
is strictly monotone increasing.  At the same time, again for all sufficiently large $u$,
 $\sigma \mapsto (\psi_{\kappa}^2 \star \Phi_{\tau/\sigma}) ( u ) $ 
 is strictly monotone decreasing in $\sigma$.
 Also, let $H_{\mu}$ denote  the CDF of a point mass at $\mu$,
then $H_{\mu/\sigma} = S^\sigma H_{\mu}$. Consequently, 
$\sigma \mapsto S^\sigma H_{\mu}$ is increasingly 
concentrated (rather than spread) as $\sigma$ increases.
It follows that, for large enough $\mu > 0$,
\[
  \sigma \mapsto  \int (\psi_{\kappa}^2 \star \Phi_{\tau/\sigma}) ( u ) \de H_{\mu}^\sigma(u) 
    = (\psi_{\kappa}^2 \star \Phi_{\tau/\sigma}) ( \mu/\sigma )
\]
is monotone decreasing in $\sigma$.
Similarly, for large enough $\mu > 0$,
\[
  \sigma \mapsto  \int (\psi_{\kappa}' \star \Phi_{\tau/\sigma}) ( u ) \de H_{\mu}^\sigma(u) 
    = (\psi_{\kappa}' \star \Phi_{\tau/\sigma}) ( \mu/\sigma )
\]
is monotone increasing in $\sigma$.

Because
\[
   A(\psi_{\kappa}, \tilde{H}^\sigma) = (\psi_{\kappa}^2 \star \Phi_{\tau/\sigma}) ( \mu/\sigma )
\]
and
\[
 B(\psi_{\kappa}, \tilde{H}^\sigma) =  (\psi_{\kappa}' \star \Phi_{\tau/\sigma}) ( \mu/\sigma ),
\]
and the decompositions $I+II$ and $III+IV$, our claims about the behavior of the RHS's in these
displays, for large $\mu$, imply the needed monotonicities of $\sigma \mapsto A(\psi_{\kappa}, \tilde{F}^\sigma)$
and $\sigma \mapsto B(\psi_{\kappa}, \tilde{F}^\sigma)$. \qed

\subsection*{Proof of Lemma \ref{lem:monotoneLambda}}
\def\cK{{\cal K}}
Putting $\bar{\sigma}(m,\kappa,\eps) \equiv \sqrt{1 + \tau_\infty^2(m,\eps,\kappa)}$,
we note that
\[
   \blambda(m,\kappa,\eps) = \kappa \cdot \bar{\sigma}(m,\kappa,\eps).
\]
Now
\[
  \bar{\sigma}(m,\kappa,\eps)^2 = (1 + \frac{\bar{V}(\bkappa,\eps)/m}{1- \bar{V}(\bkappa,\eps)/m}) =  \frac{1}{1- \bar{V}(\bkappa,\eps)/m}.
\]
By direct evaluation,
the function $\kappa \mapsto \bar{V}(\kappa,\eps)$
is at first strictly decreasing on $(0,\infty)$ to a minimum at the Huber minimax parameter
$\kappa^*(\eps)$, after which it is strictly increasing, tending to infinity
as $\kappa \goto  \infty$.

Consequently, on the interval $\kappa \in (\kappa^*(\eps),\infty)$, the function 
$\kappa \mapsto \bar{V}(\kappa,\eps)$ is strictly increasing.
On the interval $\cK_+ = (\ukappa^*(\eps),\kappa^+(\eps;m))$ 
 the function $\kappa \mapsto \bar{V}(\bkappa(\kappa),\eps)$ is likewise
 strictly increasing. Hence on $\cK_+$
 $\kappa \mapsto \bar{\sigma}$ is strictly increasing, and so also is $\kappa \cdot \bar{\sigma}$.

Fix $m_0 >  \bar{V}(0,\eps) = \frac{\pi}{2(1-\eps)} $.
For each $m > m_0$,
 $\bar{\sigma}(m,0,\eps) < \infty$, and this is the largest
that $\bar{\sigma}(m,\kappa,\eps)$ ever gets on $\kappa \in (0,\ukappa^*(\eps))$.
On the interval $\cK_- = (0, \ukappa^*(\eps))$
the function $\kappa \mapsto \bar{V}(\bkappa(\kappa),\eps)$ is bounded and has bounded derivative.
It follows that, as $m \goto \infty$
\[
  \sup_{\kappa \in \cK_-}  | \bar{\sigma}(m,\kappa,\eps)  - 1|  \goto 0, \qquad m \goto \infty;
\]
and also
\[
  \sup_{\kappa \in \cK_-}  | \frac{\partial}{\partial \kappa } \bar{\sigma}(m,\kappa,\eps)  -  0| \goto 0, \qquad m \goto \infty,
\]
together implying 
\[
  \sup_{\kappa \in \cK_-}  | \frac{\partial}{\partial \kappa } \blambda(m,\kappa,\eps)  -  1 | \goto 0, \qquad m \goto \infty,
\]
yielding $\frac{\partial}{\partial \kappa } \blambda(m,\kappa,\eps) > 0$ throughout $\cK_-$. 

We have shown that $\bar{\lambda}$ is strictly increasing, as
a function of $\kappa$, throughout the whole domain $(0,\kappa^+(m,\eps)) = \cK_- \cup \cK_+$. \qed

\subsection*{Proof of Theorem  \ref{lem:minmaxThm}}

We first remark that for any specific $\lambda > 0$,
\begin{eqnarray*}
      \cV_m^\circ(\lambda,\bar{F}_{\eps})  &=& \cV_m(\kappa(m,\lambda,\bar{F}_{\eps}), \bar{F}_{\eps}),  \\
                                           &=& m \cdot \tau_\infty^2(m,\kappa(m,\lambda,\bar{F}_{\eps}),\bar{F}_{\eps}),\\
                                          &\geq& m \cdot \min_\kappa \tau_\infty^2(m,\kappa,\bar{F}_{\eps}),\\
                                           &=& m \cdot \bbtau_\infty^2(m,\ukappa^*(\eps),\eps) \\
                                             &=&  \cV_m(\kappa^*, \bar{F}_\eps) = \cV_m^*(\eps). 
\end{eqnarray*}
and so
\[
          \inf_\lambda \sup_{F \in \cF_\eps}  \cV_m^\circ(\lambda,\bar{F}_{\eps}) \geq  \cV_m^*(\eps).
\]

We complete the argument  by showing that
\[
         \sup_{F \in \cF_\eps}  \cV_m^\circ(\lambda,\bar{F}_{\eps}) \leq  \cV_m^*(\eps),
\]
or in other words:
\[
\cV_m^\circ(\lambda^*,F)  \leq \cV_m^*(\eps), \qquad \forall F \in \cF_\eps .
\]
To show this, we need merely to show that for each $(\kappa,F)$ yielding 
an instance where $\lambda(m,\kappa,F) = \lambda^*(m,\eps,\ukappa^*(\eps))$,
we have
\beq \label{eq:tau-ineq}
\tau_\infty^2(m,\kappa,F) \leq \bbtau_\infty^2(m,\ukappa^*(\eps),\eps),
\eeq
since then
\begin{eqnarray*}
      \cV_m^\circ(\lambda^*,F)  &=& \cV_m(\kappa,F) \\
                                           &=& m \cdot \tau_\infty^2(m,\kappa,F),\\
                                            & \leq & m \cdot \bbtau_\infty^2(m,\ukappa^*(\eps),\eps) \\
                                             &=&  \cV_m(\kappa^*, \bar{F}_\eps) = \cV_m^*(\eps). 
\end{eqnarray*}
Suppose that  $\kappa \geq \ukappa(\eps)$,  then from
\begin{eqnarray*}
      \kappa \cdot \sqrt{1 + \tau_\infty^2(m,\kappa,F)}  & = & \lambda(m,\kappa,F) \\
      & = & \lambda^* \\
      & = &   \ukappa^*(\eps)  \cdot \sqrt{1 + \btau_\infty^2(m,\ukappa^*(\eps),\eps)}
\end{eqnarray*}
we conclude that
\begin{eqnarray*}
      \sqrt{1 + \tau_\infty^2(m,\kappa,F)}   \cdot \frac{\kappa}{\ukappa^*(\eps)} & = &  \sqrt{1 + \btau_\infty^2(m,\ukappa^*(\eps),\eps)}
\end{eqnarray*}
and since $ \frac{\kappa}{\ukappa^*(\eps)} \geq 1$, we indeed obtain (\ref{eq:tau-ineq}).

To finish, we argue that $\kappa < \ukappa^*(\eps)$ can never arise in a pair $(\kappa,F)$ obeying $\lambda(m,\kappa,F) = \lambda^*$.
By the monotonicity property of Lemma \ref{lem:monotoneLambda}, if we have $\kappa < \ukappa^*(\eps)$,
\begin{eqnarray*}
   \sup_{F \in \cF_\eps} \lambda(m,\kappa,F)  &=& \bar{\lambda}(m,\kappa,\eps) \\
      & < & \bar{\lambda}(m,\ukappa^*(\eps),\eps) ,
\end{eqnarray*}
proving that it can never happen that $\kappa(m,\lambda^*,F) < \ukappa^*(\eps)$, for any $F \in \cF_\eps$. \qed

\bibliographystyle{amsalpha}
\bibliography{all-bibliography}

\providecommand{\bysame}{\leavevmode\hbox to3em{\hrulefill}\thinspace}
\providecommand{\MR}{\relax\ifhmode\unskip\space\fi MR }
\providecommand{\MRhref}[2]{%
  \href{http://www.ams.org/mathscinet-getitem?mr=#1}{#2}
}
\providecommand{\href}[2]{#2}
\begin{thebibliography}{BBEKY13}

\bibitem[BBEKY13]{bean2013optimal}
Derek Bean, Peter~J Bickel, Noureddine El~Karoui, and Bin Yu, \emph{Optimal
  {M}-estimation in high-dimensional regression}, Proceedings of the National
  Academy of Sciences \textbf{110} (2013), no.~36, 14563--14568.

\bibitem[Blo74]{Bloomfield}
Peter Bloomfield, \emph{On the distribution of residuals from a fitted linear
  model}, Princeton Department of Statistics Technical Report 56, 1974.

\bibitem[BRT09]{BickelEtAl}
P.~J. Bickel, Y.~Ritov, and A.~B. Tsybakov, \emph{{Simultaneous analysis of {
  L}asso and {D}antzig selector}}, Annals of Statistics \textbf{37} (2009),
  1705--1732.

\bibitem[BvdG11]{buhlmann2011statistics}
Peter B{\"u}hlmann and Sara van~de Geer, \emph{Statistics for high-dimensional
  data}, Springer-Verlag, 2011.

\bibitem[CT07]{candes2007dantzig}
Emmanuel Candes and Terence Tao, \emph{The {D}antzig selector: Statistical
  estimation when p is much larger than n}, The Annals of Statistics (2007),
  2313--2351.

\bibitem[DH83]{donoho1983notion}
David~L Donoho and Peter~J Huber, \emph{The notion of breakdown point}, A
  Festschrift for Erich L. Lehmann (1983), 157--184.

\bibitem[DM13]{donoho2013high}
David Donoho and Andrea Montanari, \emph{{High Dimensional Robust
  {M}-Estimation: Asymptotic Variance via Approximate Message Passing}}, {\sf
  arXiv:1310.7320} (2013).

\bibitem[DM15]{andreaHuberRigor}
\bysame, \emph{{Approximate Message Passing Algorithm for the {H}uber
  {M}-estimator: Validity of State Evolution}}, {\sf unpublished}, 2015.

\bibitem[DW99]{daniel1999fitting}
Cuthbert Daniel and Fred~S Wood, \emph{Fitting equations to data: computer
  analysis of multifactor data}, John Wiley \& Sons, Inc., 1999.

\bibitem[EKBBL13]{karoui2013robust}
Noureddine El~Karoui, Derek Bean, Peter~J Bickel, and Bin Lim, Chingwayand~Yu,
  \emph{On robust regression with high-dimensional predictors}, Proceedings of
  the National Academy of Sciences \textbf{110} (2013), no.~36, 14557--14562.

\bibitem[Ham74]{hampel1974influence}
Frank~R Hampel, \emph{The influence curve and its role in robust estimation},
  Journal of the American Statistical Association \textbf{69} (1974), no.~346,
  383--393.

\bibitem[HR09]{huber2009robust}
Peter~J Huber and Elvezio~M Ronchetti, \emph{Robust statistics}, Wiley, 2009.

\bibitem[Hub64]{HuberMinimax}
P.J. Huber, \emph{Robust estimation of a location parameter}, The Annals of
  Mathematical Statistics \textbf{35} (1964), no.~1, 73--101.

\bibitem[Hub73]{huber1973robust}
Peter~J Huber, \emph{Robust regression: asymptotics, conjectures and {M}onte
  {C}arlo}, The Annals of Statistics \textbf{1} (1973), no.~5, 799--821.

\bibitem[Kar13]{karoui2013cavity}
Noureddine~El Karoui, \emph{Asymptotic behavior of unregularized and
  ridge-regularized high-dimensional robust regression estimators: rigorous
  results}, {\sf arXiv:1311.2445} (2013).

\bibitem[Por84]{portnoy1984asymptotic}
Stephen Portnoy, \emph{{Asymptotic behavior of M-estimators of $p$ regression
  parameters when $p^2/n$ is large. I. Consistency}}, The Annals of Statistics
  (1984), 1298--1309.

\bibitem[Ser10]{Serdobolski}
V.~Serdobolskii, \emph{Multivariate statistical analysis: A high-dimensional
  approach}, Kluwer Academic Publishers, 2010.

\end{thebibliography}

\end{document}